\documentclass[oneside,12pt]{article} 

\usepackage{amsmath,amssymb}
\usepackage{amsthm}  
\usepackage{mathrsfs,upref} 
\usepackage{graphicx}
\usepackage{indentfirst}

\usepackage[dvips]{hyperref}
\hypersetup{
            colorlinks=true,
            linkcolor=blue,
            anchorcolor=green,
            citecolor=red,
            bookmarksnumbered,
            bookmarksopen,
            bookmarksopenlevel=2
            }

\numberwithin{equation}{section}

\newtheorem {definition}{Definition}[section] 
\newtheorem {theorem}[definition]{Theorem}
\newtheorem {corollary}[definition]{Corollary}
\newtheorem {lemma}[definition]{Lemma}
\newtheorem {notation}[definition]{Notation}

\newtheorem {proposition}{Proposition}[section]
\newtheorem {remark}{Remark}[section]
\newtheorem {illustration}{Illustration}[section]

\newcommand{\sldef}[1]{\begin{definition}{\textsl{#1}}\end{definition}}
\newcommand{\slthm}[1]{\begin{theorem}{\textsl{#1}}\end{theorem}}
\newcommand{\slcor}[1]{\begin{corollary}{\textsl{#1}}\end{corollary}}
\newcommand{\sllem}[1]{\begin{lemma}{\textsl{#1}}\end{lemma}}
\newcommand{\slnot}[1]{\begin{notation}{\textsl{#1}}\end{notation}}

\newcommand{\slprop}[1]{\begin{proposition}{\textsl{#1}}\end{proposition}}
\newcommand{\slrem}[1]{\begin{remark}{\textsl{#1}}\end{remark}}

\newcommand{\slprf}[1]{\begin{proof}{\textsl{#1}}\end{proof}} 

\newcommand{\isomeanvalue}[3]{{\overline {{\displaystyle #1}_{\scriptstyle #2}}}\bigl|_{#3}}
\newcommand{\isomeanvalueII}[3]{{\overline {{\displaystyle #1}_{\scriptstyle #2}}}\bigl|_{#3}^{II}}  
\newcommand{\isomeanvalueIII}[3]{{\overline {{\displaystyle #1}_{\scriptstyle #2}}}\bigl|_{#3}^{III}}

\hfuzz=\maxdimen
\begin{document}

\title{\Large{\textbf{ON THE ISOMORPHIC MEANS AND\\ COMPARISON INEQUALITIES}}}
\author {Tim Y. Liu
        \thanks{
                Liu Yuan(Tim) (1975-Dec-15) \newline
                \hspace*{1.8em}Field of Research: Mathematical Analysis \newline  
                \hspace*{1.8em}Email: liu\_yuan@aliyun.com \newline
                \hspace*{1.8em}Finished Date: 2018-Jun-06 \hspace{12em} Accepted Date: \newline  
        } \\
        {\small Hengbao Co.,Ltd.}
}
\date{\small Rev. Aug 2023}

\maketitle

\begin{abstract}
Based on collection of bijections, variable and function are extended into ``isomorphic variable'' and ``dual-variable-isomorphic function'', then mean values such as arithmetic mean and mean of a function are extended to ``isomorphic means''. 7 sub-classes of isomorphic mean of a function are distinguished. Comparison problems of isomorphic means are discussed. A sub-class(class V) of isomorphic mean of a function related to Cauchy mean value is utilized for generation of bivariate means e.g. quasi-Stolarsky means. Demonstrated as an example of math related to ``isomorphic frames'', this paper attempts to unify current common means into a better extended family of means.
\vspace{0.2cm}

\vspace{0.2cm}

MSC2010 Subject Classification: 26A24 
\end{abstract}

\section{Introduction}
The quasi-arithmetic mean \cite{BULLENPS} or generalized $f$-mean $f^{-1}\big(\frac{1}{n}\sum_{i=1}^n f(x_i)\big)$ is a generalization of simple means (of numbers) using a function $f$. Typical special cases are arithmetic mean, geometric mean and power mean etc.

As a byproduct of the study of a generalization of convex function in article \cite{LIUY} in the early 2000s, the author claims having independently discovered an analogous ``2-D'' version of so-called ``generalized $g,h$-mean of a function $F$ defined on $[a,b]$'' instead(denoted by $M_F$), which typically can be in the formula:
\begin{equation} \label{equ:IsoMeanofFunIntro}
    M_F = h^{-1}\Big(\frac{1}{g(b)-g(a)}\int_{g(a)}^{g(b)}h\bigl(F(g^{-1}(u))\bigl)\mathrm{d}u\Big),
\end{equation}
where $g,h$, comparable to above $f$, are 2 bijections each applied on 2 variables of $F$ respectively(i.e. its independent variable and dependent variable, henceforth called ``dual-variable-'' throughout this work).

With its essential intermediate value property (IVP)(Theorem \ref{thm:IVPofIsoMean}) applicable to an ordinary function $F$, this complex form of mean is uncovered as a HUGE class of mean of a function. It is defined as ``isomorphic mean of a function'' (Definition \ref{defin:IsoMeanOfFun}) and 7 sub-classes of it are distinguished in Section \ref{subsec:IsoMeanClass}, such its special cases and derivations are abundant and they cross with many existing concepts of mean values. This focused ``2-D'' version mean of a function along with the ``1-D'' generalized-$f$ mean are unified to a sweeping class of ``Isomorphic Means''(as in the title) due to the same nature explained below.

In Algebra, two algebraic systems based on 2 sets are mapped by a bijection with which two operations in one system each have the same structures, such they are ``isomorphic'' and these 2 operations are ``images'' to each other. The author borrows the concepts and terms of ``isomorphic'' and ``image'' in naming, renaming and explaining of the above 2 MA objects and many others in an unified way.
\begin{itemize}\setlength{\itemsep}{-0.1em}
  \item[1).] We introduce the ``isomorphic number'' as a new concept which structure is mapped from that of a variable by a bijection $f$. Then generalized $f$-mean is the image of the mean of several isomorphic number instances by the inverted bijection. Thus it is redefined as ``isomorphic mean (of numbers) generated by $f$''.
  \item[2).] Extending to the scope of 2 variables(dual-variable), we introduce the ``dual-variable-isomorphic(DVI) function of $f$'' generated by a pair of bijections, which structure is ``co-mapped'' by the pair from that of the original $f$ (of 1 independent variable). With regards to $f\colon D\to M$, it has the form of
    \begin{equation}
        \varphi\colon=(h\circ f\circ g^{-1})\colon E\to N ~~~(E=g(D),~N=h(M)).
    \end{equation}
Then (\ref{equ:IsoMeanofFunIntro}) is rightly the image of the mean of the DVI function of $F$ therein, thus it is named as ``isomorphic mean of $F$ generated by $g,~h$''.
    \item[3).] As a matter of fact, the ``isomorphic mean of a function'' can be directly derived from ``isomorphic mean of numbers", as in Section \ref{subsec:IsoMeanFunEquaDev}.
    \item[4).] Naturally these structures can be further generalized to $n$-dimensions for functions of $(n-1)$ variables. See the very last formula (\ref{equ:IsoMeanValueN}).
\end{itemize}

For efficient discussions and evolutions, the related $n$ monotone bijections as a single mapping needed by definitions is refined as a basic concept of ``isomorphic frame'', such the isomorphic means elaborated in this paper as a theory of generalization is just an example of mathematics related to isomorphic frames.

Besides the basic properties and subtle classifications of isomorphic means, this paper features the following findings and results:
\begin{itemize}\setlength{\itemsep}{-0.1em}
  \item[1).]  The abundant special cases of isomorphic means introduced in the paper. Among others are ``elastic mean(E)'' derived from Economics, and the unfamous ``geometric mean of a function(G)'' as the sibling of the former. Because isomorphic means of a function may derive mean values of numbers, (e.g. when $g,h$ are power functions and $f(x)=x$) there obtained a class of ``quasi-Stolarsky means'' of the following form:
      \begin{equation} \label{equ:QuasiStoIntro}
        Q_{p,q}(a,b)= \biggl( \frac{p(b^{p+q}-a^{p+q})}{(p+q)(b^p-a^p)}\biggl)^{1/q},
      \end{equation}
      in the sense that this mother form yields all major children forms(see Section \ref{sec:GenBivarMean}) of Stolarsky means(\cite{YANGZHCvCompHM}, pp629).
  \item[2).]  This paper discusses the relations and differentiations of isomorphic means to(from) the outstanding ``Cauchy mean values'' and their conversions(Section \ref{sec:RelationIsoMeanCauchy}). The major finding is that the IVP of isomorphic means works with ordinary functions, while the IVP of Cauchy mean values applies to derivable functions; and the Cauchy mean values can only be confidently converted to and from the so-called class V of isomorphic means of derivable generators. Especially for mean of a function, isomorphic mean is a concept of better origin and perspective, better identification and classification, more coverage and more natural generalization.
  \item[3).]  The comparison problems of isomorphic means are roughly solved, which take up a major part of this paper. The main results are Theorem \ref{thm:IsomeanC5CompNess} and Theorem \ref{thm:IsomeanC5CompSuff} derived from Losonczi's theorems in \cite{LOSONLCAUCHYCOMP}, for class V comparison scenarios; and the Theorem \ref{thm:IsoMeanC2CompCvx}, which is for comparison of delicate class II and features very original proof based on the very properties of isomorphic mean itself; and along with those theorems derived by help of monotonicity and convexity conditions for specific comparison scenarios. Typical examples are the comparisons between the geometric mean(G) and the elastic mean(E) of a function, among many others. There is also a Lemma \ref{lem:IsoMeanComp} quoted from article \cite{LIUY}, and Theorem \ref{thm:IsoMeanCompNew}, which all delicately compares 2 isomorphic means of numbers(generalized-$f$ means).
\end{itemize}

In the general sense, the above-mentioned ``isomorphic mean of a function'' has the formal name of ``the dual-variable-isomorphic(DVI) mean of a function'', as defined in Definition \ref{defin:IsoMeanOfFun}.

\pagebreak

\section{Basics and preliminaries}
There are 3 basic concepts to be introduced. Namely ``Isomorphic frame'', ``Isomorphic number and Isomorphic variable'', and ``Dual-variable-isomorphic function''.
\subsection{Isomorphic frame}
\subsubsection{Definition}
\sllem{\label{lem:OrderSetbeBJ}Let $X_1,...,X_n,U_1,...,U_n\subseteq\mathbb{R}$, ~$X=X_1\times ...\times X_n$, $~U=U_1\times ...\times U_n$ and $g_i\colon X_i\to U_i(i=1,...,n)$ be $n$ bijections. The ordered set $\{g_1,...,g_n\}$ that can map $\forall x=(x_1,...,x_n)\in X$ to $u=(g_1(x_1),...,g_n(x_n))\in U$ is a bijection, if it is treated as a function with $X$ being its domain and $U$ being the range(the image).}

\slprf{Firstly $\forall x=(x_1,...,x_n)\in X, ~y=(y_1,...,y_n)\in X$ satisfying $(g_1(x_1),...,g_n(x_n))=(g_1(y_1),...,g_n(y_n))$, then due to the $n$ bijections, $x_1=y_1,~...,~x_n=y_n \Rightarrow x=y$(injective). ~Secondly $\forall u=(u_1,...,u_n)\in U$, there always be an $x=(g_1^{-1}(u_1),...,g_n^{-1}(u_n))\in X$ such that $\{g_1,...,g_n\}(x)=u$(surjective).}

\sldef{The ordered set $\{g_1,...,g_n\}$ given by Lemma \ref{lem:OrderSetbeBJ} as a bijection as well as a collection, denoted by $\mathscr{I}\{g_1,...,g_n\}$, ~is called an $n$-dimensional isomorphic frame. It is further expressed as the following notational forms:
\begin{equation} 
\begin{split}
   \mathscr{I}\{g_1,...,g_n\}:X\to U &= [X~\sharp~U]_{g_1,...,g_n}\\
                                     &= [X_1\times ...\times X_n~\sharp~ U_1\times...\times U_n]_{g_1,...,g_n}\\
                                     &= [X_1, ...,X_n~\sharp~ U_1,...,U_n]_{g_1,...,g_n}.
\end{split}
\end{equation}
$X$ is called the base frame of the isomorphic frame(or ``the base'' for short), and $~U$ is called the image frame of the isomorphic frame(or ``the image''). The bijection $g_i(i=1,...,n)$ is called a (the $i$th) dimensional mapping(DM) of $\mathscr{I}\{g_1,...,g_n\}$.
}
\slnot{We also write $\mathscr{I}^{-1}\{g_1,...,g_n\}$ for the inversion of $\mathscr{I}\{g_1,...,g_n\}$.}
\slthm{$\mathscr{I}\{g_1^{-1},...,g_n^{-1}\}=(\mathscr{I}\{g_1,...,g_n\})^{-1}(=\mathscr{I}^{-1}\{g_1,...,g_n\})$.}

(Proof omitted.)

\subsubsection{Embedding in isomorphic frame}
\slnot{Let $\mathscr{I}\{g_1,...,g_n\} = [X_1\times ...\times X_n~\sharp~U_1\times...\times U_n]_{g_1,...,g_n}$ and $D\subseteq X_1\times ...\times X_n$, $E\subseteq U_1\times ...\times U_n$. \\
\indent (i) $D$ is said to be embedded in the base of $\mathscr{I}\{g_1,...,g_n\}$, for which way we use the sign ``$\vee$'' to write $D\vee\mathscr{I}\{g_1,...,g_n\}$, or $D\vee\big(\mathscr{I}\{g_1,...,g_n\}=[X_1\times ...\times X_n~\sharp~U_1\times...\times U_n]_{g_1,...,g_n}\big)$, etc.\\
 \indent (ii) $E$ is said to be embedded in the image of $\mathscr{I}\{g_1,...,g_n\}$, and for this we write $E\vee\mathscr{I}\{g_1^{-1},...,g_n^{-1}\}$.
}

\slthm{$\mathscr{I}\{g_1,...,g_n\}(D)\vee\mathscr{I}\{g_1^{-1},...,g_n^{-1}\}$ if $D\vee\mathscr{I}\{g_1,...,g_n\}$.
}
This is because the image of $D$ is a subset of the image of the isomorphic frame.

\subsubsection{Bonding on isomorphic frame}
\slnot{Let $\mathscr{I}\{g_1,...,g_n\} = [X_1\times ...\times X_n~\sharp~U_1\times...\times U_n]_{g_1,...,g_n}$.\\
\indent (i) If function $f:D\to M$ of $(n-1)$ variables is such that $D\subseteq X_1\times ...\times X_{n-1}$ ($D\vee\mathscr{I}\{g_1,...,g_{n-1}\}$) and the range $M\subseteq X_n$ ($M\vee\mathscr{I}\{g_n\}$), then $f$ is said to be bonded on the base of $\mathscr{I}\{g_1,...,g_n\}$. For this we use the sign ``$\wedge$'' to write $(f:D\to M) \wedge\mathscr{I}\{g_1,...,g_n\}$, or the alike. \\
\indent (ii) If $D\subseteq U_1\times ...\times U_{n-1}$ and $M\subseteq U_n$, then $f$ is said to be bonded on the image of $\mathscr{I}\{g_1,...,g_n\}$, and for this we write $(f:D\to M) \wedge\mathscr{I}\{g_1^{-1},...,g_n^{-1}\}$, $f \wedge\mathscr{I}^{-1}\{g_1,...,g_n\}$ or the alike.}

\slrem{In the extreme case $n=1$, the above notation applies to a single variable, e.g. $x\in M\subseteq X$, or $x\in M\subseteq U$, which are treated as functions of 0 variables bonded on (the base or the image of) the isomorphic frame.}

\slnot{Let $\mathscr{I}\{g\}= [X~\sharp~U]_g$ be 1 dimensional, and $k$-tuple $\underline{x}\in X$, then $\underline{x}$ is said to be either embedded in or bonded on the base of $\mathscr{I}\{g\}$, written as e.g. $\underline{x}\vee[X~\sharp~U]_g$, or $\underline{x}\wedge[X~\sharp~U]_g$.
}

\subsubsection{About embedding and bonding}
Embedding and bonding are 2 basic styles the objects of MA ``attach or tie'' to isomorphic frames. A set embedded in an isomorphic frame is a part of the latter while a function bonded on an isomorphic frame is not that way strictly. Most definitions in the paper will be based on these 2 concepts.

From now on, all the isomorphic frames are with $\underline{strictly ~monotone~ bijections}$, i.e. strictly monotone (invertible) real functions, as their dimensional mappings unless otherwise specified. These isomorphic frames are denoted as $\mathscr{I}_m\{g_1,...,g_n\}$.

\subsection{Isomorphic number and isomorphic variable}
\sldef{\label{defin:IsoNumVar} With 1 dimensional $\mathscr{I}_m\{g\} = [X~\sharp~U]_g$, $\forall x\in X$, $u=g(x)\in U$ is called the isomorphic number of $x$ generated by mapping $g$(or by $\mathscr{I}_m\{g\}$); In terms of variables, $u$ is called the isomorphic variable of $x$ generated by mapping $g$. $u$ is specially denoted as $u=\varphi(x:g)$ or $u=\varphi(x:\mathscr{I}_m\{g\})$.}

\slrem{Any real number $x$ is an isomorphic number of itself generated by identity mapping ($y=x$).}
An example of isomorphic number in Physics is Conductance G with regards to Resistance R, since $G=1/R$.

With Definition \ref{defin:IsoNumVar} $u$ is already called the ``isomorphic number of $x$'' without an operation defined, because with such $u$ it is ready and easy to introduce 2 genuine operations with ``isomorphism''. For examples:
\begin{itemize}\setlength{\itemsep}{-0.1em}
\item[1).] If there is an operation of arithmetic ``+'' in $U$, that $\forall u_1,u_2\in U$, $\exists u_1+u_2\in U$, then we can define a binary operation denoted by $[\ .+.\ ]_g$ on $X$, such that $\forall x_1,x_2\in X$, $\exists x_3\in X$ which holds $x_3=[x_1+x_2]_g=g^{-1}(g(x_1)+g(x_2))$.
\item[2).] Similarly, if there's an operation on $U$ that computes the mean of $u_1$, $u_2$, there is also a mapped operation on $X$, which computes a ``special mean'': $\bar x$ of $x_1,x_2$. 
\end{itemize}
More generalized, $u$ and $x$ are ``2 conjugate variables embedded in $[X~\sharp~U]_g$'' or ``2 imaging functions of 0 variables bonded on $[X~\sharp~U]_g$''. 

\subsection{Dual-variable-isomorphic function}
With a function of 1 variable bonded on the base of a 2-D isomorphic frame, on the image of the latter bonded is the so called ``dual-variable-isomorphic function''.

\subsubsection{Definition}

\sldef{\label{defin:IsoFun} Let $(f\colon D\to M)\wedge\big(\mathscr{I}_m\{g,h\} = [X,Y~\sharp~U,V]_{g,h}\big)$ and $E=g(D),~N=h(M)$,
\begin{itemize}
\item Function  $(h\circ f\circ g^{-1})\colon E\to N$ is defined to be the dual-variable-isomorphic(DVI) function of $f$ generated by mapping $g,h$(or generated by $\mathscr{I}_m\{g,h\}$). It is denoted by $\varphi(f:g,h)$, or $\varphi(f:\mathscr{I}_m\{g,h\})$,
\begin{equation} 
    \varphi(f:g,h)\colon=(h\circ f\circ g^{-1})\colon E\to N.
\end{equation}
\item $g$ and $h$ are called the independent variable's generator mapping(function) or dimensional mapping(IVDM), and dependent variable's generator mapping(function) or dimensional mapping(PVDM) of $\varphi$ respectively.
\item $E\subseteq U$ is called the isomorphic domain of $D$ generated by mapping $g$; $N\subseteq V$ is called the isomorphic range of $M$ generated by mapping $h$.
\end{itemize}
}

\begin{figure}[htb]
\begin{center}
\includegraphics[keepaspectratio=true,scale=0.4]{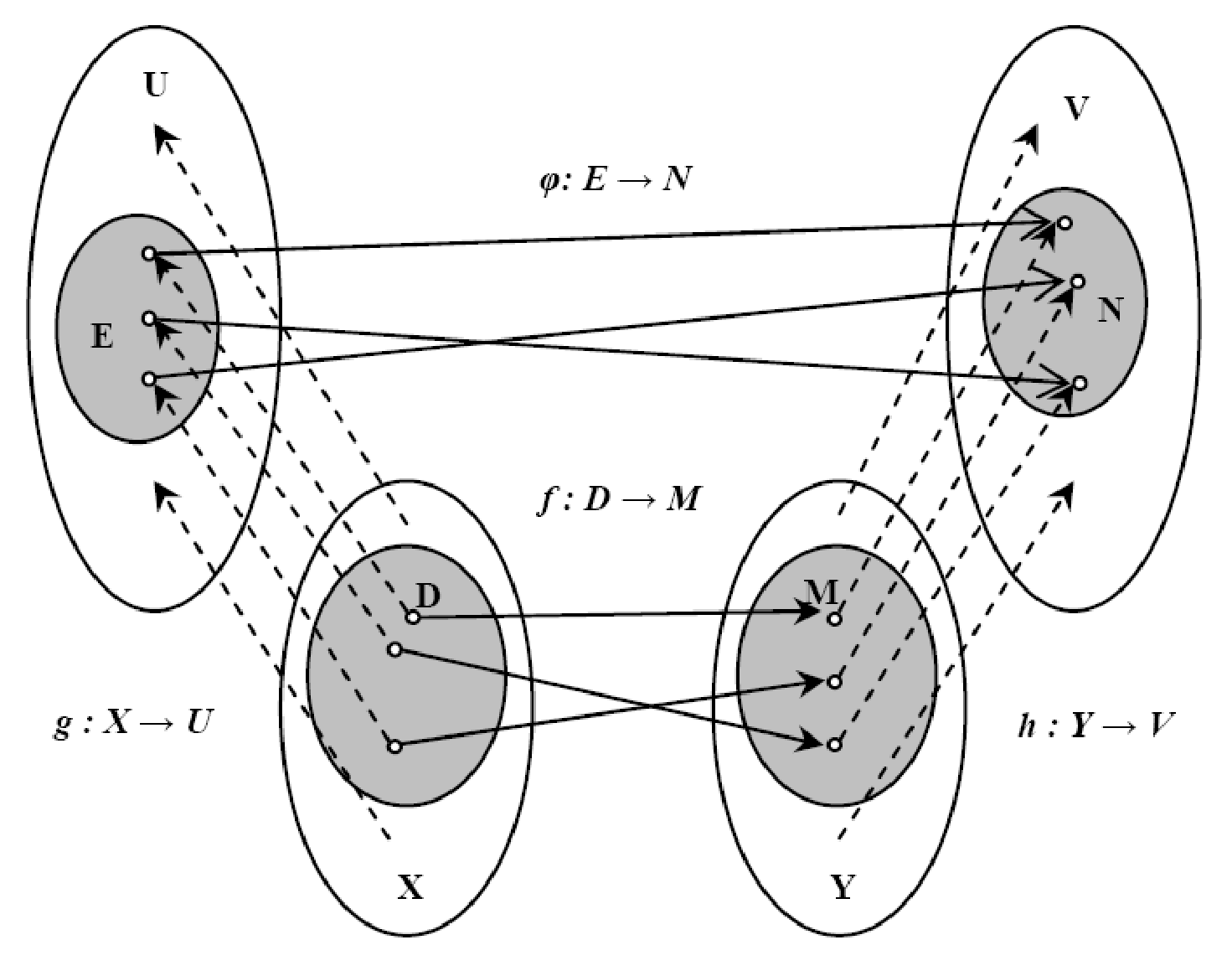}
\end{center}
\begin{illustration}{
\begin{center}
\label{illus:IsoFun}Function $f$ and its DVI function $\varphi\colon E\to N$ \\bonded on a 2-D isomorphic frame
\end{center}
}\end{illustration}
\end{figure}

\slthm{$\varphi(f:g,h) \wedge \big(\mathscr{I}_m\{g^{-1},h^{-1}\} = [U,V~\sharp~X,Y]_{g^{-1},h^{-1}}\big)$.}
That is, $\varphi(f:g,h)$ is bonded on the image of $\mathscr{I}_m\{g,h\}$. Proof omitted.

Illustration \ref{illus:IsoFun} is the visual impression of the definition and the theorem, where $f\wedge\mathscr{I}_m\{g,h\}$ and $\varphi\wedge\mathscr{I}_m\{g^{-1},h^{-1}\}$. A hidden mapping between $f$ and $\varphi$ is formed by $\mathscr{I}_m\{g,h\}$ under which these 2 functions are deemed as either (i) 2 similar elements in different spaces, (ii) 2 sets of relations of same cardinal number, or (iii) 2 unary operations of same structure but each spanning across a pair of number sets.

In traditional terms, the range $M$ here is the image of $f$, and $Y$ is the codomain. The following expression sometimes may also be used for DVI functions, with which the range $N\subseteq V$ is implied, where $V$ is also the codomain.
\begin{equation} 
    \varphi(f:g,h)\colon=(h\circ f\circ g^{-1})\colon E\to V.
\end{equation}

\subsubsection{Sub-classing of DVI function}\label{subsubsec:IsoFunSpecialCase}
\sldef{In view of Definition \ref{defin:IsoFun}, the following special cases can be considered:
\begin{itemize}\setlength{\itemsep}{-0.1em}
  \item[1).] Let $g$ be identity mapping, $\varphi\colon=(h\circ f)\colon D\to N$ is called the dependent-variable-isomorphic(PVI) function of $f$ generated by $h$.
  \item[2).] Let $h$ be identity mapping, $\varphi\colon=( f\circ g^{-1})\colon E\to M$ is called the independent-variable-isomorphic(IVI) function of $f$ generated by $g$.
  \item[3).] Let $Y=X$, $h=g$, $\varphi\colon=(g\circ f\circ g^{-1})\colon E\to N$ is called the same-mapping dual-variable-isomorphic function of $f$ generated by $g$.
  \item[4).] In general case, $h \neq g$, $\varphi\colon=(h\circ f\circ g^{-1})\colon E\to N$ is called the (general) dual-variable-isomorphic function of $f$ generated by $g,h$.
  \item[5).] Let $f(x)=x$, $\varphi\colon=(h\circ g^{-1})\colon E\to N$ is called the dual-variable-isomorphic function of identity generated by $g,h$.
  \item[6).] Let $g$, $h$ both be identity mappings, $\varphi\colon=(h\circ f\circ g^{-1})\colon E\to N$  is equivalent to $f\colon D\to M$. Both domains are $D$, and both ranges are $M$.  i.e. $f$ is a special dual-variable-isomorphic function of itself generated by identities.
  \item[7).] For monotone function $f\colon D\to M(range~M)$, $f^{-1}$ is the dual-variable-isomorphic function of $f$ generated by $f$, $f^{-1}$. 
  \\Above special cases also could be treated as 7 sub-classes of DVI function.
\end{itemize}
}

\subsubsection{Anti-dual-variable-isomorphic function}

\sldef{Let $\varphi\colon=(h\circ f\circ g^{-1})\colon E\to N$  be the DVI function of $f\colon D\to M$. Then $f$ is called the anti-dual-variable-isomorphic function of $\varphi$.}

\slrem{With respect to its DVI function $\varphi$ generated by mapping $g$, $h$, \,$f\colon D\to M$ can be represented by $f=(h^{-1}\circ \varphi \circ g) :D\to M$. Obviously $f$ is the DVI function of $\varphi$ generated by $g^{-1}$, $h^{-1}$. }

This can be observed in Illustration \ref{illus:IsoFun}.

The following are 3 useful special cases of DVI functions.

\subsubsection{V-scaleshift: Vertical scale and shift of a function}
\sldef{For a real $f\colon D\to M$ and 2 constants $k\ne0$, $C$, define the function
\begin{equation}
    V_{ss}\big(f:k,C\big)\colon=~~v=kf(x)+C ~~\big(v\in k(M)+C\big)
\end{equation}
a V-scaleshift of $f$ with scale $k$ and shift $C$. Define the set
\begin{equation}
    \mathbb{V}f = \{ V_{ss}\big(f:k,C\big):k,C\in \mathbb{R},k\ne0\}
\end{equation}
the V-scaleshift space of $f$.
}
\sllem{1).$f\in \mathbb{V}f$; ~2).$\mathbb{V}g=\mathbb{V}f$, if $g\in \mathbb{V}f$.}
\slrem{A V-scaleshift of $y=f(x)$ is a special case of dependent-variable-isomorphic function of $f$, where the PVDM is $v=ky+C$.}
\sllem{Suppose $f\colon D\to M$ be (strictly) convex or (strictly) concave on $D$, then $V_{ss}\big(f:k,C\big)$ and $f$ are of the same (strict) convexity  if $k>0$ or of the opposite (strict) convexity if $k<0$.}

Proof is omitted.

\slnot{In this paper, $\mathbb{V}x$($\mathbb{V}y$) is denoting the V-scaleshift space of identity mapping $g(x)=x$($h(y)=y$).}

\subsubsection{H-scaleshift: Horizontal scale and shift of a function}
\sldef{For a real $f\colon D\to M$ and 2 constants $k\ne0$, $C$, define the function
\begin{equation}
    H_{ss}\big(f:k,C\big)\colon=~~y=f\big(\frac1{k}(u-C)\big) ~~\big(u\in k(D)+C\big)
\end{equation}
an H-scaleshift of $f$ with scale $k$ and shift $C$. Define the set
\begin{equation}
    \mathbb{H}f = \{ H_{ss}\big(f:k,C\big):k,C\in \mathbb{R},k\ne0\}
\end{equation}
the H-scaleshift space of $f$.
}
\sllem{1).$f\in \mathbb{H}f$; ~2).$\mathbb{H}g=\mathbb{H}f$, if $g\in \mathbb{H}f$.}
\slrem{An H-scaleshift of $y=f(x)$ is a special case of independent-variable-isomorphic function of $f$, where the IVDM is $u=kx+C$.}
\sllem{\label{lem:HssSameConvex}Suppose $f\colon D\to M$ be (strictly) convex or (strictly) concave on $D$, then $H_{ss}\big(f:k,C\big)$ is of the same (strict) convexity on $k(D)+C$ as $f$ on $D$.
}
\slprf{Suppose $f$ is convex, then $\forall u_1,u_2\in k(D)+C$, $\forall \lambda\in [0,1]$ $\exists x_1=(u_1-C)/k, x_2=(u_2-C)/k\in D$, such that $f(\lambda x_1+(1-\lambda)x_2)\leq \lambda f(x_1)+(1-\lambda )f(x_2)$. This $\Rightarrow$ $f(\lambda(u_1-C)/k+(1-\lambda)(u_2-C)/k)=f((\lambda u_1+(1-\lambda)u_2)/k-C)\leq \lambda f((u_1-C)/k)+(1-\lambda )f((u_2-C)/k)$, which means $H_{ss}\big(f:k,C\big)$ is also convex. While $f$ has other (strict) convexities, $H_{ss}\big(f:k,C\big)$ will also copy.}

\sllem{\label{lem:HssVssMirror}Let $g,h$ be both invertible. 1).If $g\in \mathbb{H}h$, then $g^{-1}\in \mathbb{V}(h^{-1})$. 2).If $g\in \mathbb{V}h$, then $g^{-1}\in \mathbb{H}(h^{-1})$.}
Proof is omitted.

\subsubsection{HV-scaleshift: Horizontal and vertical scale and shift of a function}
\sldef{For a real $y=f(x):D\to M$ and constants $p\ne0$, $Q$, $k\ne0$, $L$,  define the function
\begin{equation}
\begin{split}
    HV_{ss}&\big(f:p,Q;k,L\big)\colon=~~v=kf\big(\frac1{p}(u-Q)\big)+L \\
    &\big(u\in p(D)+Q,~v\in k(M)+L\big)
\end{split}
\end{equation}
an HV-scaleshift of $f$ with scale $p,k$ and shift $Q,L$. Define the set
\begin{equation}
    \mathbb{HV}f = \{ HV_{ss}\big(f:p,Q;k,L\big):p,Q,k,L\in \mathbb{R},p\ne0,k\ne0\}
\end{equation}
the HV-scaleshift space of $f$.
}
\sllem{1).$f\in \mathbb{HV}f$; ~2).$\mathbb{HV}g=\mathbb{HV}f$, if $g\in \mathbb{HV}f$.}
\slrem{
An HV-scaleshift of $y=f(x)$ is a special case of DVI function of $f$:
\begin{equation}
    HV_{ss}\big(f:p,Q;k,L\big)=\varphi(f:g,h),
\end{equation}
where the DMs $g,h$ are defined by $~u=g(x)=px+Q, ~v=h(y)=ky+L$.
}

\section{Isomorphic weighted mean of numbers}\label{sec:IsoWghtMean}
After preparation, the generalized $f$-mean will be re-defined as the ``isomorphic (weighted) mean'' with isomorphic frame involved.

\subsection{Isomorphic weighted mean}
\sldef{Let $\mathscr{I}_m\{g\}=[X~\sharp~U]_g$ be 1 dimensional, $n$-tuple($n\geq2$) $\underline{x}\wedge[X~\sharp~U]_g$ and $\underline{u}$ being their respective isomorphic numbers. With positive $\underline{p}$ satisfying $\sum_{i=1}^np_i=1$ if  $\sum_{i=1}^np_iu_i\in U$, then $g^{-1}\bigl(\sum_{i=1}^np_iu_i\bigl)$ is called the isomorphic weighted mean of $n$-tuple (numbers) $\underline x$ generated by $g$(or by $\mathscr{I}_m\{g\}$). Here it is denoted by $\overline{x,p_R}|_g$ (or $\overline{x_{\{i\}},p_R}|_g$, or simpler $\overline{x,p}|_g$),
\begin{equation}\label{equ:IsoWgtMean}
    \overline{x,p_R}|_g=g^{-1}\bigl(\sum_{i=1}^np_i g(x_i)\bigl).
\end{equation}
$g$ is called the generator mapping, or dimensional mapping of the isomorphic weighted mean.
}

\slrem{The subscript $R$ of $p$ indicates $\underline p$ are the ``relative'' weights (fractions always add up to 1), as comparing to another type of weights known as Frequency Numbers, for which case the definition and formula will be slightly different as in \cite{BULLENPS}. For simplicity, in this paper we use $\overline{x,p}|_g$ which agrees series $\underline p$ are relative (fractional) weights.
}
\sldef{Especially when $p_1=\ldots =p_n=1/n$, the isomorphic weighted mean is called the isomorphic mean of $n$-tuple $\underline x$ generated by $g$, it is denoted by $\overline{x_{\{i\}}}|_g$ or $\overline{x}|_g$,
\begin{equation}\label{equ:IsoMean}
    \overline{x}|_g=g^{-1}\bigl(\frac 1n\sum_{i=1}^ng(x_i)\bigl).
\end{equation}
}

\slrem{In simple words, isomorphic (weighted) mean is the inverse image of the (weighted) mean of $n$ number of $\varphi(x_i:g)$.}

Among many already known properties, the following are some key properties of the mean.
\subsubsection{Property of mean value}
\slthm{\label{thm:XiofIsoWgtMean} With $n$-tuple $\,\underline{x}\wedge(\mathscr{I}_m\{g\}=[X~\sharp~U]_g)$ that are not all equal, where $g$ is continuous on interval $X$, and with weights $~\underline{p}$, there is an unique $\xi\in (\min\{\underline x\}, \max\{\underline x\})\subseteq X$ such that
\begin{equation}
    \xi=\overline{x,p}|_g=g^{-1}\bigl(\sum_{i=1}^np_i g(x_i)\bigl).
\end{equation}
}

This theorem reflects a basic property of above defined isomorphic weighted mean. However in the definition of quasi-arithmetic mean in pp266 of \cite{BULLENPS}, ~$g$ being continuous is prerequisite thus the existence of the mean is ensured. The proof is omitted.

\subsubsection{Property of monotonicity}
\slthm{\label{thm:MonotofIsoWgtMean} Any participating number($x_i$)'s value increasing will result in the isomorphic weighted mean's increasing.}

Thus any one's decreasing results in the mean's decreasing. This property is obvious with (\ref{equ:IsoWgtMean}) as $g$ and $g^{-1}$ are always of the same strict monotonicity.

\subsubsection{Invariant value with vertical scale and shift of DM}
\slthm{$\overline{x,p}|_h=\overline{x,p}|_g$, for  $h\in \mathbb{V}g$.}
\slprf{$\overline{x,p}|_h=g^{-1}\biggl(\Bigl(-C+\sum_{i=1}^np_i \big(kg(x_i)+C\big)\Bigl)/k\biggl)=\overline{x,p}|_g$.}

\subsection{Comparison of 2 isomorphic (weighted) means of numbers}
For same $n$-tuple($n\geq2$) $\underline{x}$ and weights $\underline{p}$, 2 sets of methods are proposed to compare these 2 means:
\begin{equation}\label{equ:IsoMean2}
\begin{split}
    \overline{x,p}|_g=g^{-1}\bigl(\sum_{i=1}^np_i g(x_i)\bigl),
    ~~~\overline{x,p}|_h=h^{-1}\bigl(\sum_{i=1}^np_i h(x_i)\bigl). \nonumber
\end{split}
\end{equation}

Note: In this paper regarding the convexity of a function, ``convex'' means convex to lower, and ``concave'' means convex to upper.

\subsubsection{Comparison method derived from Jensen's inequality}
The following theorem may already be available, though we reproduce it here with our own simple proof, also serving to the purpose of easier cross-reference.
\linespread{1.2}
\slthm{\label{thm:IsoMeanCompOld} Suppose with $n$-tuple $~\underline{x}\wedge [X~\sharp~U]_g$, $[X~\sharp~V]_h$, and weights $\underline p$, there exists $\overline{x,p}|_g$, $\overline{x,p}|_h$. Let $D=[\min\{\underline{y}\}, ~\max\{\underline{y}\}]$ where tuple $\underline{y}=h(\underline{x})$. Then
\\\indent1). $\overline{x,p}|_g \geq \overline{x,p}|_h$, if $g$ is increasing on $X$ and $g(h^{-1})$ is convex on $D$;
\\\indent2). $\overline{x,p}|_g \leq \overline{x,p}|_h$, if $g$ is increasing on $X$ and $g(h^{-1})$ is concave on $D$;
\\\indent3). $\overline{x,p}|_g \leq \overline{x,p}|_h$, if $g$ is decreasing on $X$ and $g(h^{-1})$ is convex on $D$;
\\\indent4). $\overline{x,p}|_g \geq \overline{x,p}|_h$, if $g$ is decreasing on $X$ and $g(h^{-1})$ is concave on $D$.
\\For all cases, the equality holds only if $n$-tuple $\underline{x}$ are all equal.
}
\linespread{1.0}
\slprf{Assume $g$ is strictly increasing, the inequality between $\overline{x,p}|_g$ and $\overline{x,p}|_h$ is equivalent to the inequality between $\sum_{i=1}^np_ig(h^{-1}(y_i))$ and $g(h^{-1}(\sum_{i=1}^np_i y_i))$. According to Jensen's inequality \cite{KOSMALA}, if $g(h^{-1})$ is a convex function(convex to lower) on $[\min\{\underline{y}\}, ~\max\{\underline{y}\}]$ then
\begin{equation}
\sum_{i=1}^np_ig(h^{-1}(y_i)) \geq g(h^{-1}(\sum_{i=1}^np_iy_i)).
\end{equation}
At the same time, $\overline{x,p}|_g \geq \overline{x,p}|_h$. While $g$ may also be decreasing which reverses above-mentioned equivalence of inequalities, or $g(h^{-1})$ may be concave which reverses above Jensen's inequality, there totals 4 cases by the combinations, which summarize as the theorem. As for equality, it corresponds to how Jensen's inequality behaves similarly.}

For instance of Theorem \ref{thm:IsoMeanCompOld}, let $g(x)=\sin x$, $~y=h(x)=\ln x$, $x_i\in (0, \pi/2]$. $g$ is strictly increasing, $g(h^{-1}(y))=sin(e^y)$. Simple calculation concludes that when $1/x\geq \tan x$, $\mathrm{d}^2g(h^{-1}(y))/\mathrm{d}y^2\geq0$, which indicates $g(h^{-1}(y))$ is convex on $D$. Solved the inequality with
iterative method we get $x\in (0,0.8603...)$.

Then according to Theorem \ref{thm:IsoMeanCompOld} case 1), $\overline{x,p}|_{\sin x}\geq \overline{x,p}|_{\ln x}$ when $x_i\in(0,0.8603...)$. And contrarily with Theorem \ref{thm:IsoMeanCompOld} case 2), $\overline{x,p}|_{\sin x}\leq \overline{x,p}|_{\ln x}$ when $x_i\in(0.8603..., \pi/2)$.

\subsubsection{Comparison method of differential criteria}
\linespread{1.2}
\slthm{\label{thm:IsoMeanCompNew}With $n$-tuple $~\underline{x}\wedge [X~\sharp~U]_g$, $~\underline{x}\wedge [X~\sharp~V]_h$, and weights $\underline p$, where $X$ is an interval on which $g,h$ are monotone and derivable, $g'\ne 0, h'\ne 0$,
\\\indent1). $\overline{x,p}|_g\geq \overline{x,p}|_h$, if $|g'/h'|$ is increasing on $X$;
\\\indent2). $\overline{x,p}|_g\leq \overline{x,p}|_h$, if $|g'/h'|$ is decreasing on $X$.\\
For both cases, the equality holds only if $n$-tuple $\underline x$ are all equal.}
\linespread{1.0}
\slprf{Let $y=h(x)$, then $g(h^{-1}(y))' = (g'/h')(x)\circ h^{-1}(y)$. With Theorem \ref{thm:IsoMeanCompOld} to distinguish 8 cases: \\
\indent1). In the case $g'>0$ thus $g$ is increasing, and $(g'/h')(x)$ is increasing with $h'>0$ thus $g(h^{-1}(y))$ is convex, ~which implies $|g'/h'|=|g'|/|h'|=g'/h'$ is increasing, then $\overline{x,p}|_g \geq \overline{x,p}|_h$.\\
\indent2). In the case $g'>0$ thus $g$ is increasing, and $(g'/h')(x)$ is increasing with $h'<0$ thus $g(h^{-1}(y))$ is concave, which implies $|g'/h'|=|g'|/|h'|=g'/(-h')$ is decreasing, then $\overline{x,p}|_g \leq \overline{x,p}|_h$.\\
\indent3). In the case $g'<0$ thus $g$ is decreasing, and $(g'/h')(x)$ is increasing with $h'>0$ thus $g(h^{-1}(y))$ is convex, which implies $|g'/h'|=|g'|/|h'|=(-g')/h'$ is decreasing, then $\overline{x,p}|_g \leq \overline{x,p}|_h$.\\
\indent4). through 8). ...(omitted) \\
All the omitted 5 cases are analogous and 8 cases are finally merged without overlapping or confliction into these 2 cases of the theorem.}

\subsubsection{The nature of the comparisons from a new perspective} \label{subsubsec:natureCompIsoMean}
Also based on the DVI function, article \cite{LIUY} is yet another theory on extension of convex function. Its Theorem 2 is about how to use the monotonicity of a derivative-like function $(h\circ f)'/g'$ to compare an isomorphic-weighted-mean-involved inequality, so as to identify the extended convexity of $f$. A special case (where $f(x)=x$)  proposed as its Corollary 1, is quoted as below Lemma \ref{lem:IsoMeanComp} (some letters are changed):
\linespread{1.2}
\sllem{\label{lem:IsoMeanComp}Let $g, h$ be strictly monotone continuous and derivable functions on an open interval $X$, and $g'\ne0, ~h'\ne0$. $h'/g'$ is monotone. There are $n$-tuple($n\geq2$)$~\underline x\in X$, and $n$-tuple positive $\underline p$ adding up to 1. Among functions $h'/g', h, g$,
\\\indent1). if odd numbers (1 or 3) are increasing then $\overline{x,p}|_g\leq \overline{x,p}|_h$;
\\\indent2). if odd numbers (1 or 3) are decreasing then $\overline{x,p}|_g\geq \overline{x,p}|_h$.
\\For both cases, the inequalities are equal only if $n$-tuple $\underline x$ are all equal.
}
\linespread{1.0}
It's easy to see Theorem \ref{thm:IsoMeanCompNew} is just an improved version of the Lemma \ref{lem:IsoMeanComp}. In article \cite{LIUY}, with simple steps Lemma \ref{lem:IsoMeanComp} facilitates the proof of the power mean inequality:

\emph{For positive $n$-tuple($n\geq2$) $\underline x$ and real number $p>q$,
\begin{equation}\label{equ:PowerMean}
    \overline{x}|_{x^p}\geq \overline{x}|_{x^q},
\end{equation}
the inequality is equal only when $n$-tuple $\underline x$ are all equal.}

Another example by Lemma \ref{lem:IsoMeanComp} in pp85 of article \cite{LIUY} is: $g(x)=\sinh(x), ~h(x)=\cosh(x)$. On $\mathbb{R}^+$ $g,h$ are increasing, and $h'/g'=(e ^x-e^{-x})/(e^x+e^{-x})=1-2/(e^{2x}+1)$ is increasing, thus
\begin{equation}
    \overline{x,p}|_{\sinh x}\leq \overline{x,p}|_{\cosh x}~(x_i>0).
\end{equation}

 On the other hand, in Analysis the method of using the monotonicity of $f'$ (to compare an inequality) to determine the convexity of $f$ is also a special case of Theorem 2 in \cite{LIUY}, since convexity is actually the basic case of the ``extended convexity'' in article \cite{LIUY}. And the Lemma \ref{lem:IsoMeanComp} is a ``cousin'' method of the former. (Similarly in \cite{LIUY} is a third cousin method for ``geometrically convexity''.)

Correspondingly $h'/g'$ in Lemma \ref{lem:IsoMeanComp} is a special form of $(h\circ f)'/g'$ when $f$ is identity, as ``cousin'' to the derivative of $y=x$ which is also a special form of $(h\circ f)'/g'$ when $h,g$ are identities. But being not always 1 the former does have monotonicity indeed.

\textbf{In this sense, the inequality between 2 isomorphic means generated by $g$ and by $h$ is merely the indication of the ``extended convexity generated by $g,h$'' of $y=x$.}

\subsubsection{An impression of comparisons}
The author has used the methods to compare isomorphic means pairs among 5 common generator mappings: $y=\sin x$, $~y=1/x$, $~y=\ln x$, $~y=x$, and $y=e^x$ for $x_i\in[-\pi/2,\pi/2]$. A graphical view of the results is as below:

\begin{figure}[htb]
\begin{center}
\includegraphics[height=80mm,width=140mm]{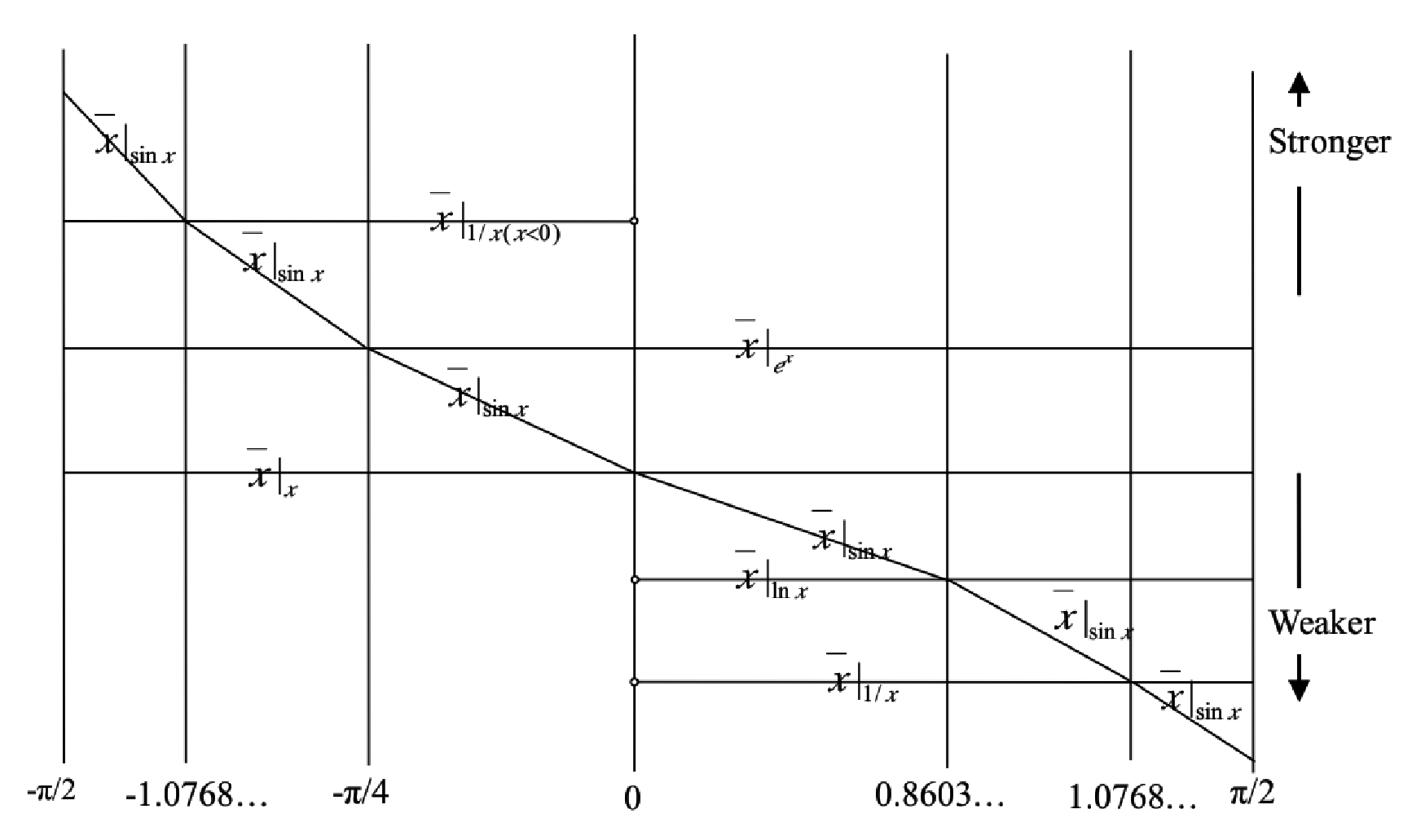}
\end{center}
\begin{illustration}{
\begin{center}
\label{illus:IsoMean5Comp}An impression of comparison of 5 types of isomorphic means
\end{center}
}\end{illustration}
\end{figure}

Taking the 2nd vertical partition from left for example, it shows that for $x_i\in(-1.0768..., -\pi/4)$, $\overline x|_{1/x}\geq \overline x|_{\sin x}\geq\overline x|_{e^x}\geq\overline x|_{x}$. An imaginary falling slope is drawn for $\overline x|_{\sin x}$, which displays a rightward trend of ``being weaker'', reflecting an interesting property of the sine function. While other inequalities such as arithmetic mean being stronger than geometric mean can also be observed.\\

In the next section, one can see that the redefinition of ``generalized $f$-mean'' to isomorphic mean of numbers is in harmony with ``isomorphic mean of a function'', and the former is also deemed as a basis of the latter.

For this reason, here we give an informal classification to the isomorphic (weighted) mean of numbers: ``the isomorphic mean class 0''.(as more classes of isomorphic mean will come up next.)

\section{Isomorphic mean of a function}

When the mean values of a function bonded on an isomorphic frame is concerned, the dual-variable-isomorphic(DVI) mean of a function can be similarly introduced.

\subsection{The definition of DVI mean of a function}
\sldef{\label{defin:IsoMeanOfFun} Let $f:D\to M$ be bounded on interval $D$ and $f\wedge\big(\mathscr{I}_m\{g,h\}=[X,Y~\sharp~U,V]_{g,h}\big)$. $g$ is continuous on $D$, interval $E=g(D)$, and $h$ is continuous on $[\inf M, \sup M]$. If there exists $M_\varphi \in V$, being the mean value of $\varphi\colon= (h\circ f \circ g^{-1})\colon E\to V$, then $h^{-1}(M_\varphi)\in Y$ is called the (dual-variable-) isomorphic mean(DVI mean) of $f$ on $D$ generated by $g, h$ (or generated by $\mathscr{I}_m\{g,h\}$), denoted by $\isomeanvalue{f}{D}{g,h}$, $\isomeanvalue{f}{\,}{g,h}$ or $M_f|_{g,h}$.\\
\begin{equation}\label{equ:IsoMeanValue0}
\begin{split}
    \isomeanvalue{f}{D}{g,h} &=
    h^{-1}\biggl(\frac{\int_{\scriptscriptstyle E}h\big(f(g^{-1}(u))\big)\mathrm{d}u}
    {\int_{\scriptscriptstyle E}\mathrm{d}u}\biggl). \\
    \biggl(~~&=h^{-1}\bigg(\frac{\int_{\scriptscriptstyle E}h\big(f(g^{-1}(u))\big)\mathrm{d}(-u)}
    {\int_{\scriptscriptstyle E}\mathrm{d}(-u)}\bigg).~~\biggl)
\end{split}
\end{equation}
It's stipulated that $\forall a\in D$, $\isomeanvalue{f}{[a,a]}{g,h}=f(a)$. ~$g,h$ are called the independent variable's generator mapping or dimensional mapping(IVDM), and the dependent variable's generator mapping or dimensional mapping(PVDM) of the isomorphic mean respectively.
}
\slrem{If $M_\varphi$ exists, then there must be $M_\varphi\in h([\inf M, \sup M])\subseteq V$, because $h([\inf M, \sup M])=[\inf N, \sup N]$ ($N=h(M)$) in which $M_\varphi$ must resides.} 

\slrem{In simple words, isomorphic mean of a function is the inverse image of the mean of $\varphi(f:g,h)$.}
Intended to be the extension of mean of a function on an interval, this article defines the isomorphic mean of an ``ordinary'' function rather on an interval, than on a general real number set. Another ad-hoc requirement in addition to the ``bonding'' is $[\inf M, \sup M]\subseteq Y$.

\slnot{Taking $[a, b]$ as the form of $D$, the formula (\ref{equ:IsoMeanValue0}) turns into
\begin{equation}\label{equ:IsoMeanValue}
\begin{split}
    \isomeanvalue{f}{[a,b]}{g,h} &= h^{-1}\biggl(\frac{1}{g(b)-g(a)}\int_{g(a)}^{g(b)}h\bigl(f(g^{-1}(u))\bigl)\mathrm{d}u\biggl),\\
    \biggl(~~ &= h^{-1}\bigg(\frac{1}{g(a)-g(b)}\int_{g(b)}^{g(a)}h\bigl(f(g^{-1}(u))\bigl)\mathrm{d}u\bigg) ~~\biggl)
\end{split}
\end{equation}
with which whether $g(b)>g(a)$ is disregarded since the result of the formula is the same with that having $g(a)$ and $g(b)$ exchanged their positions. And if interval $D$ has infinite endpoint(s), the limit form of (\ref{equ:IsoMeanValue}) may be considered, with $a$ and/or $b$ approaching infinity.
}

Whereas if $f$ is not bounded, which falls out of scope of Definition \ref{defin:IsoMeanOfFun}, as long as the value of (\ref{equ:IsoMeanValue0}) exists, it could be considered the generalized isomorphic mean of an unbounded function.

\subsection{\label{subsec:ExiCompatIsoMean}The existence of \texorpdfstring{$M_f|_{g,h}$}{Mf|g,h} and tolerance of the definition}
With Definition \ref{defin:IsoMeanOfFun}, the existence of $M_f|_{g,h}$ depends on the existence of the mean of $\varphi(f:g,h)$. The following theorem is sufficient but not necessary.
\slthm{\label{thm:IsomeanExist1}The $\isomeanvalue{f}{D}{g,h}$ exists with any applicable $g,~h$, if $f$ is continuous on a close interval $D$.
}
\slprf{Such $f,g,h$ imply a continuous $\varphi(f:g,h)$ on a close interval $E$ with a convex range $h(M)$. This means the numerator integral in the formula (\ref{equ:IsoMeanValue0}) exists and the dominator is non-zero real, such $M_\varphi$ exists in $h(M)\subseteq V$. These lead to an unique $\isomeanvalue{f}{D}{g,h}$ in $M\subseteq Y$.
}

While one expects for $f$ being continuous, the range $M$ being convex(an interval) and $M_f|_{g,h} \in M$, among others the definition however allows for
\begin{itemize} \setlength{\itemsep}{-0.2em}
  \item[(1).] $f$ is not continuous with jump discontinuity, but $M_f|_{g,h}$ exists. e.g.:
   \begin{equation} 
    \begin{split}
        f(x)&=\left\{
        \begin{aligned}
        1, ~~(x&\in [0,1]) \\
        3. ~~(x&\in (1,2]) \\
        \end{aligned} \right.\\
        g(x)&=2x, ~h(y)=y+2.
    \end{split}  
   \end{equation}
   With above, $M_f|_{g,h}=2$.
  \item [(2).] $f$ is not continuous with essential discontinuity, but $M_f|_{g,h}$ exists. e.g.
  If $f$ is a discontinuous but bounded Darboux function, especially being derivative of another continuous function $F$, and being Riemann integrable, then with the same $g,h$ as of above, $M_f|_{g,h}$ exists.
  \item [(3).] $M$ is not convex, and $M_f|_{g,h}\notin M$; but $M_f|_{g,h}\in Y$ as with above case 1).
  \item [(4).] $M$ is not convex, but $M_f|_{g,h}\in M\subseteq Y$.
  \item [(5).] The definition is not contradict with: $f$ is not continuous, $M_f|_{g,h}$ does not exist. e.g. $f(x)=\pi, x\in[0,2]$ and $x$ is irrational; $f(x)=3.14, x\in[0,2]$ and $x$ is rational, with the same $g,h$ as of above.
\end{itemize}
While it does not allow for $M_f|_{g,h}\in Y$ and $M_f|_{g,h}\notin [\inf M, \sup M]$, which will be proved later.

Another sufficient condition for the existence of isomorphic mean is:
\slthm{\label{thm:IsomeanExist2}The $\isomeanvalue{f}{D}{g,h}$ exists if $f$ is Riemann integrable on a close interval $D=[a,b]$ with an applicable $\mathscr{I}_m\{g,h\}$ such $g$ satisfies one of the following:
\begin{itemize} \setlength{\itemsep}{-0.2em}
  \item[1).] $g$ is derivable on $(a,b)$;
  \item[2).] $g$ is a convex or concave function on $[a,b]$;
  \item[3).] $g$ is absolutely continuous on any close sub-interval of $(a,b)$.
\end{itemize}
}
\slprf{Any of the 3 additional conditions adding to the strictly monotone and continuous $g$ ensures $f\circ g^{-1}$ is also integrable on $E$. With $h$ being continuous, $h\circ f\circ g^{-1}$ is Riemann integrable. Meanwhile the dominator is non-zero real, hence $M_\varphi$ exists in $h([\inf M, \sup M])\subseteq V$. This leads to an unique $\isomeanvalue{f}{D}{g,h}$ in $[\inf M, \sup M]\subseteq Y$.
}

Regarding the Riemann integrability of composite functions, the discussions can be founded in \cite{LUJ} and \cite{HUANGQL} which are applicable to e.g. $h\circ f\circ g^{-1}$ here. According to \cite{HUANGQL} any of additional conditions 1) 2) applying to $g$ will lead to condition 3) that will map the set of discontinuity points of $f$ on $[a,b]$ of Lebesgue measure 0 to a counterpart of $f\circ g^{-1}$ on $E$ of Lebesgue measure 0.

With above we claim that isomorphic means of a function are generally available in 2-D isomorphic frames with derivable DMs, for ORDINARY functions on close intervals, which (i) are Riemann integrable; (ii) may be somewhere discontinuous; (iii) have not to be monotone.

\subsection{An equivalent derivation of DVI mean of a function} \label{subsec:IsoMeanFunEquaDev}
Isomorphic mean of a function can also be derived through the limit of the isomorphic weighted mean of numbers expressed in integral form. In Definition \ref{defin:IsoMeanOfFun} assuming $D$ is a close interval $[a,~b]$, then due to monotone and continuous $g$, $E$ is a close interval $[g(a),~g(b)]$(disregarding whether $g(b)>g(a)$). Let $\tau=f\circ g^{-1}$. Do a $n$-tuple partitions of $E$ similar to is done with a definite integral. With each partition $i~(i\in{1,2,...,n}$), there is a corresponding value $\tau(\xi_i)$ for its tagged point $\xi_i$. Then the ratio of $\Delta u_i$($=u_i-u_{i-1}$, where $u_i$ is the end point of each partition towards $g(b)$, $u_{i-1}$ is the other end point of the same partition) over $g(b)-g(a)$, is taken as $\tau(\xi_i)$'s weight $w_i$. We compute the isomorphic weighted mean of $n$-tuple $\tau(\xi_i)$ generated by $h$, denoted by $M_{\tau(\xi_i)}$,
\begin{equation} 
\begin{split}
    M_{\tau(\xi_i)}  &= h^{-1}\biggl( \sum_{i=1}^{n}{\frac{u_i-u_{i-1}}{g(b)-g(a)}h(\tau(\xi_i))}\biggl) \label{equ:PartialSum}\\
        &= h^{-1}\biggl( \frac{1}{g(b)-g(a)}\sum_{i=1}^{n}{h(\tau(\xi_i))\Delta u_i}\biggl).
\end{split}  
\end{equation}

In order for $\tau(\xi_i)$ to enumerate all possible values of $\tau$ on $E$(such $f$ enumerates all on $[a,b]$), let $\|\Delta\|(=\max\{\Delta u_i\})\to0$, i.e. the partitions become infinitely thin, such
\begin{equation} 
\begin{split}
    \lim _{\|\Delta\|\to0}M_{\tau(\xi_i)} &=  h^{-1}\biggl(\frac{1}{g(b)-g(a)}\int_{g(a)}^{g(b)}h\big(\tau(u)\big)\mathrm{d}u\biggl)  \\
    &=  h^{-1}\biggl(\frac{1}{g(b)-g(a)}\int_{g(a)}^{g(b)}h\big(f(g^{-1}(u))\big)\mathrm{d}u\biggl). \label{equ:PartialSumLimit}
\end{split}
\end{equation}
This limit value is deemed as ``function $\tau$'s isomorphic weighted mean generated by mapping $h$''. It is just the same in form as (\ref{equ:IsoMeanValue}), that is function $f$'s dual-variable-isomorphic mean generated by $g,h$.

A special case of above is, the partitions are of all equal size, which is $\frac{1}{n}|g(b)-g(a)|$, such the weights are all $1/n$. Then (\ref{equ:PartialSumLimit}) can be deemed as an integral evolution of isomorphic mean of numbers, which has the same value if exists.

If $E$ is infinite or $D$ is open or half open, the limit forms of (\ref{equ:PartialSumLimit}) will be included cases of (\ref{equ:IsoMeanValue0}).

\subsection{Basic properties of DVI mean of a function}

\subsubsection{Property of intermediate value(IVP)}
\slthm{\label{thm:IVPofIsoMean}For a bounded $f\colon D\to M$($M$ is range) with its isomorphic mean $\isomeanvalue{f}{D}{g,h}$ in the applicable isomorphic frame $\mathscr{I}_m\{g,h\}$ , it holds
\begin{equation} \label{Inequ:IsoMeanFunGenIneq}
    \inf M \leq \isomeanvalue{f}{D}{g,h} \leq \sup M;
\end{equation}
especially if exists $\min\{M\}= \inf M$, $\max\{M\}=\sup M$, then
\begin{equation} \label{Inequ:IsoMeanFunGenIneq2}
   \min\{M\} \leq \isomeanvalue{f}{D}{g,h} \leq \max\{M\}.
\end{equation}
}
\slprf{Let $N=h(M)$. (i) if $M_\varphi\in N \subseteq V$, then due to the bijection $h$, $h^{-1}(M_\varphi)\in M \subseteq[\inf M, ~\sup M]$; (ii) if in general case $M_\varphi \in N$ is not to be considered, with (\ref{equ:PartialSum}), $h$ being continuous on $[\inf M, \sup M]$ by definition, and Theorem \ref{thm:XiofIsoWgtMean}, it holds
\begin{equation}
    \inf M \leq \min\{\tau(\xi_i)\}\leq M_{\tau(\xi_i)} \leq \max\{\tau(\xi_i)\}\leq \sup M
\end{equation}
while $\|\Delta\|\to0$. Thus $\inf M \leq \isomeanvalue{f}{D}{g,h} \leq \sup M$. For all cases, if there exists $\min\{M\}, \max\{M\}$, it holds $\min\{M\} \leq \isomeanvalue{f}{D}{g,h} \leq \max\{M\}$.}
Generally this Intermediate Value Property(IVP) or mean value property of an isomorphic mean does not require the corresponding intermediate value(s) within the domain of the function.

This IVP holds with ordinary functions. It could be fundamentally accepted as an attribute of ordinary functions on the background of isomorphic frames of continuous DMs, which will make isomorphic means have more coverage of extended unique means of a function than those concepts of IVP only holding with continuous or derivable functions.

The transforming nature of the background as a single bijection also makes the isomorphic mean an extended mean value of more straightforward and genuine origin, as directly from the mean value of an transformed function. As a fact, it naturally refines the development of mean values from simpler arithmetic mean into complicated ones, so that the isomorphic mean is a natural generalization of the simple mean.

In view of these, the isomorphic mean of a function is rather a differentiated concept of IVP than a repeated one with some others. (Also see Section \ref{sec:DifferenceIsoCauchy}.)

\subsubsection{Property of monotonicity}
\slthm{\label{thm:IsoMeanMonotone}Let $f,m$ be defined on interval $D$, interval $D_2\subseteq D$. If
\begin{equation} \label{equ:MonotoneSubIneq}
m(x)\geq f(x) ~~\forall x\in D_2,
\end{equation}
and $m(x)=f(x)$ $\forall x\notin D_2$, then for any applicable $\mathscr{I}_m\{g,h\}$, it holds
\begin{equation} \label{equ:IsoMeanMonotone}
   M_m|_{g,h} \geq M_f|_{g,h}.
\end{equation}
If (\ref{equ:MonotoneSubIneq}) is strict, then (\ref{equ:IsoMeanMonotone}) is strict.
}
\slprf{While there are several cases according to the ways $D_2$ is located in $D$, we first assume both be close intervals, such $D$ can be expressed as $[a,b]$~ $(a<b)$. And arrange $D_2$ such $D_2=[a,c]$, $(a<c\leq b)$, i.e. $D_2$ locates to the left of $D$.  Let $\tau=f\circ g^{-1}$, $\mu=m\circ g^{-1}$, with the same way in Section \ref{subsec:IsoMeanFunEquaDev}, do a set of $n$-tuple partitions and tagged $\xi_i$ for $D_2=[a,c]$, another set of $k$-tuple partitions and tagged $\eta_j$ for $(c,b]$. Let
\begin{equation} 
\begin{split}
    M_{\tau(\xi,\eta)}  &= h^{-1}\biggl( \sum_{i=1}^{n}\frac{u_i-u_{i-1}}{g(b)-g(a)}h(\tau(\xi_i))+
    \sum_{j=1}^{k}\frac{v_j-v_{j-1}}{g(b)-g(a)}h(\tau(\eta_j))\biggl) \\
    M_{\mu(\xi,\eta)}  &= h^{-1}\biggl( \sum_{i=1}^{n}\frac{u_i-u_{i-1}}{g(b)-g(a)}h(\mu(\xi_i))+
    \sum_{j=1}^{k}\frac{v_j-v_{j-1}}{g(b)-g(a)}h(\mu(\eta_j))\biggl)
\end{split}
\end{equation}
\indent Above each putting in $h^{-1}$ are 2 partial sums. For convenience let's denote them $S_1$, $S_2$ for those inside $M_{\tau(\xi,\eta)}$ resp., and $T_1$, $T_2$ for those inside $M_{\mu(\xi,\eta)}$ resp., such
\begin{equation} \label{equ:IsomeanMonotone}
\begin{split}
    M_{\tau(\xi,\eta)}  = h^{-1}\big(S_1+S_2\big), ~~
    M_{\mu(\xi,\eta)}  = h^{-1}\big(T_1+T_2\big)
\end{split}
\end{equation}
\indent In the case $c=b$, $D_2=D$, $S_2,T_2$ do not exist. As $\forall x\notin D_2(i.e. ~x\in (c,b])$, $m(x)=f(x)$, such $S_2 \equiv T_2$. Meanwhile $\forall x\in D_2$, $m(x)\geq f(x)$, this makes $T_1\geq S_1 \Rightarrow T_1+T_2 \geq S_1+S_2$ with increasing $h$, or $T_1\leq S_1 \Rightarrow T_1+T_2 \leq S_1+S_2$ with decreasing $h$. In both cases $M_{\mu(\xi,\eta)} \geq M_{\tau(\xi,\eta)}$. When the partitions become infinitely thin, the inequality holds, thus it holds $M_m|_{g,h} \geq M_f|_{g,h}$.\\
 \indent For cases where a close $D_2$ is located other ways in a close $D$, the proof is analogous; Then at most we have 3 partial sums each in (\ref{equ:IsomeanMonotone}).\\
 \indent For at least one of $D_2$,$D$ and interval(s) of $D-D_2$ being open/half open intervals including above, the proofs always need to threat the endpoints carefully not to let tagged points be the endpoints that do not belongs, as fortunately finite points exclusion of tagged points does not affected the integral calculations;\\
 \indent For the cases $D_2$ and/or $D$ being infinite intervals, the proofs shall further consider the holding inequality on a selected partial finite interval, and let it hold under a limit of the endpoint(s). \\
 \indent Finally if (\ref{equ:MonotoneSubIneq}) is strict, then the corresponding inequalities within the proof including the result are strict.
}

\slcor{\label{thm:IsoMeanMonotoneCor}Let $f,m$ be defined on interval $D$, there are finite disjoint intervals $D_1,...,D_n\subseteq D$. If
\begin{equation} \label{equ:MonotoneSubIneqCor}
m(x)\geq f(x) ~~\forall x\in D_1\cup ... \cup D_n,
\end{equation}
and $m(x)=f(x)$ $\forall x\notin D_1\cup ... \cup D_n$, then for any applicable $\mathscr{I}_m\{g,h\}$, it holds
\begin{equation} \label{equ:IsoMeanMonotoneCor}
   M_m|_{g,h} \geq M_f|_{g,h}.
\end{equation}
If (\ref{equ:MonotoneSubIneqCor}) is strict, then (\ref{equ:IsoMeanMonotoneCor}) is strict.
}
This can be proved via Theorem \ref{thm:IsoMeanMonotone} and construction of $(n-1)$ bridging functions and passing the same inequality down total $(n+1)$ isomorphic means: $M_m|_{g,h} \geq...\geq M_f|_{g,h}$.

\subsubsection{Property of symmetry with endpoints of interval}
\slprop{(As implied by the definition,) The value of an isomorphic mean of a function on a close interval $[a,b]$ is invariant with $a,b$ exchanging their values in the formulae (\ref{equ:IsoMeanValue}).}

As a result, any function of the pair $(a,b)$ derived by (\ref{equ:IsoMeanValue}), maybe of different forms but of the equivalence of (\ref{equ:IsoMeanValue}), is symmetrical with $a,b$.
\slnot{For an existing $\isomeanvalue{f}{[a,b]}{g,h}$, also denote $\isomeanvalue{f}{[b,a]}{g,h}$ for the same isomorphic mean.
}
Therefore future for such $[a,b]$, $b\ge a$ is no longer a compulsory requirement.

\subsubsection{Invariant value with vertical scale and shift of dimensional mappings}
\paragraph{Invariant value with V-scaleshift of IVDM}
\slthm{\label{thm:IsoMeanInvar1} $M_f|_{m,h}=M_f|_{g,h}$ for $m\in \mathbb{V}g$.}
\slprf{
Let $u=m(x)$, $v=g(x)$, then $m^{-1}(u)=g^{-1}(v)$, ~$u=kv+C$, taking $[a,b]$ as $D$, thus
\begin{equation} \nonumber
\begin{split}
  M_f|_{m,h} &=h^{-1}\biggl(\frac{1}{kg(b)+C-kg(a)-C}\int_{kg(a)+C}^{kg(b)+C}h\big(f(m^{-1}(u))\big)\mathrm{d}u\biggl) \\
  &= h^{-1}\biggl(\frac{1}{kg(b)-kg(a)}\int_{g(a)}^{g(b)}h\big(f(g^{-1}(v))\big)\mathrm{d}(kv+C)\biggl) ~~(v=\frac{u-C}{k})\\
  &= h^{-1}\biggl(\frac{1}{g(b)-g(a)}\int_{g(a)}^{g(b)}h\big(f(g^{-1}(v))\big)\mathrm{d}v\biggl) = M_f|_{g,h}.
\end{split}
\end{equation}
While cases with other forms of $D$ are treated as holding limit forms of above.
}

\paragraph{Invariant value with V-scaleshift of PVDM}
\slthm{\label{thm:IsoMeanInvar2} $M_f|_{g,l}=M_f|_{g,h}$ for $l\in \mathbb{V}h$.}
\slprf{Taking $[a,b]$ as $D$, then
\begin{equation} \nonumber 
\begin{split}
  M_f|_{g,l} &= h^{-1}\biggl(\biggl(\bigg(\frac{1}{g(b)-g(a)}\int_{g(a)}^{g(b)}\big(kh\big(f(g^{-1}(u))\big)+C\big)\mathrm{d}u\bigg)-C\biggl)/k\biggl) \\
  &= h^{-1}\biggl(\biggl(\bigg(\frac{1}{g(b)-g(a)}\int_{g(a)}^{g(b)}kh\big(f(g^{-1}(u))\big)\mathrm{d}u+C\bigg)-C\biggl)/k\biggl) \\
  &= h^{-1}\biggl(\frac{1}{g(b)-g(a)}\int_{g(a)}^{g(b)}h\big(f(g^{-1}(u))\big)\mathrm{d}u\biggl)
  = M_f|_{g,h}.
\end{split}
\end{equation}
While cases with other forms of $D$ are treated as holding limit forms of above.
}

\paragraph{Invariant value with V-scaleshifts of both DMs}
\slcor{\label{cor:IsoMeanInvar12} $M_f|_{m,l}=M_f|_{g,h}$ for $m\in \mathbb{V}g$ and $l\in \mathbb{V}h$.}
It's a result of two-step process by previous 2 theorems.

\slnot{Also denote $M_f|_{\mathbb{V}g,\mathbb{V}h}$ for $M_f|_{g,h}$, $M_f|_{\mathbb{V}g}^{II}$ for $M_f|_g^{II}$, etc., i.e. for any IVDM or PVDM $g$, $\mathbb{V}g\Leftrightarrow g$ in denoting of isomorphic means of a function.}

\subsection{Sub-classing of DVI mean of a function} \label{subsec:IsoMeanClass}
There are 7(seven) typical sub-classes of isomorphic mean of a function, in correspondence with 7 special cases of DVI function:
\sldef{\label{Def:IsoMeanOfFunSub}With Definition \ref{defin:IsoMeanOfFun}, consider the following cases due to special $g,h,f$: 
\begin{enumerate}
  \item[1).] Let $g$ be identity, then $\varphi\colon=(h\circ f)\colon D\to V$. The isomorphic mean of $f$ is called the dependent-variable-isomorphic mean(PVI mean) of $f$ on $D$ generated by mapping $h$, or the isomorphic mean class I of $f$ on $[a, b]$ generated by $h$. It is denoted by $\isomeanvalue{f}{D}{h}$, or $M_f|_h$,
        \begin{equation}
            \isomeanvalue{f}{D}{h}
            = h^{-1}\Biggl(\frac{\int_{\scriptscriptstyle D}h\big(f(x)\big)\mathrm{d}x}
            {\int_{\scriptscriptstyle D}\mathrm{d}x}\Biggl)~(=\isomeanvalue{f}{D}{\mathbb{V}x,\mathbb{V}h}).
        \end{equation}
  If taking $[a, b]$ as $D$,
        \begin{equation}
            \isomeanvalue{f}{[a,b]}{h}
            =h^{-1}\biggl(\frac{1}{b-a}\int_a^bh\big(f(x)\big)\mathrm{d}x\biggl).
        \end{equation}
  \item[2).] Let $h$ be identity, then $\varphi\colon=(f\circ g^{-1})\colon E\to M$. The isomorphic mean of $f$ is called the independent-variable-isomorphic mean(IVI mean) of $f$ on $D$ generated by mapping $g$, or the isomorphic mean class II of $f$ on $D$ generated by $g$. It is denoted by $\isomeanvalueII{f}{D}{g}$, or $M_f|_g^{II}$,
        \begin{equation}
             \isomeanvalueII{f}{D}{g}
            = \frac{\int_{\scriptscriptstyle E}f(g^{-1}(u))\mathrm{d}u}
            {\int_{\scriptscriptstyle E}\mathrm{d}u}~(=\isomeanvalue{f}{D}{\mathbb{V}g,\mathbb{V}y}).
        \end{equation}
     If taking $[a, b]$ as $D$,
        \begin{equation}
            \isomeanvalueII{f}{[a,b]}{g}
            =\frac{1}{g(b)-g(a)}\int_{g(a)}^{g(b)}f\big(g^{-1}(u)\big)\mathrm{d}u.
        \end{equation}
     If $g$ is differentiable,
        \begin{equation}
            \isomeanvalueII{f}{[a,b]}{g}
            =\frac{1}{g(b)-g(a)}\int_a^bf(x)\mathrm{d}g(x).
        \end{equation}
  \item[3).] Let $Y=X, ~h=g$, then $\varphi\colon=(g\circ f\circ g^{-1})\colon E\to V$. The isomorphic mean of $f$ is called the same-mapping (dual-variable-)isomorphic mean(SDVI mean) of $f$ on $D$ generated by mapping $g$, or the isomorphic mean class III of $f$ on $D$ generated by $g$. It is denoted by $\isomeanvalueIII{f}{D}{g}$, $M_f|_{g,g}$, or $M_f|_{g}^{III}$,
      \begin{equation}
            \isomeanvalueIII{f}{D}{g}=
            g^{-1}\Biggl(\frac{\int_{\scriptscriptstyle E}g\big(f(g^{-1}(u))\big)\mathrm{d}u}
            {\int_{\scriptscriptstyle E}\mathrm{d}u}\Biggl)~(=\isomeanvalue{f}{D}{\mathbb{V}g,\mathbb{V}g}).
      \end{equation}
      If taking $[a, b]$ as $D$,
      \begin{equation}
            \isomeanvalueIII{f}{[a,b]}{g}
            =g^{-1}\biggl(\frac{1}{g(b)-g(a)}\int_{g(a)}^{g(b)}
            g\bigl(f(g^{-1}(u))\bigl)\mathrm{d}u\biggl).
      \end{equation}
  \item[4).] If in general case $h\ne g$, then $\varphi\colon=(h\circ f \circ g^{-1})\colon E\to V$. The isomorphic mean $\isomeanvalue{f}{D}{g,h}(=\isomeanvalue{f}{D}{\mathbb{V}g,\mathbb{V}h})$ is called (dual-variable-)isomorphic mean(DVI mean) of $f$ on $D$ generated by mapping $g, h$, or the isomorphic mean class IV of $f$ on $D$ generated by $g, h$, which formula is as (\ref{equ:IsoMeanValue0}), or as (\ref{equ:IsoMeanValue}) if taking $[a, b]$ as $D$.
  \item[5).] Let $f$ be identity, then $\varphi\colon=(h \circ g^{-1})\colon E\to V$. The isomorphic mean simplifies to the mean of one variable $x$. In this paper it is called the (dual-variable-)isomorphic mean of identity on $D$ generated by mapping $g, h$, or the isomorphic mean class V on $D$ generated by $g, h$. It is denoted by $\isomeanvalue{x}{D}{g,h}$, or $M_x|_{g,h}$,
      \begin{equation}
            \isomeanvalue{x}{D}{g,h}=
            h^{-1}\Biggl(\frac{\int_{\scriptscriptstyle E}h\big(f(g^{-1}(u))\big)\mathrm{d}u}
            {\int_{\scriptscriptstyle E}\mathrm{d}u}\Biggl)~(=\isomeanvalue{x}{D}{\mathbb{V}g,\mathbb{V}h}).
      \end{equation}
      If taking $[a, b]$ as $D$,
        \begin{equation} \label{equ:IsoMeanT5}
            \isomeanvalue{x}{[a,b]}{g,h}
            =h^{-1}\Biggl(\frac{1}{g(b)-g(a)}\int_{g(a)}^{g(b)}h\bigl(g^{-1}(u)\bigl)\mathrm{d}u\Biggl).
        \end{equation}
      If $g$ is differentiable,
        \begin{equation} \label{equ:IsoMeanT5dif}
            \isomeanvalue{x}{[a,b]}{g,h}
            =h^{-1}\Biggl(\frac{1}{g(b)-g(a)}\int_a^bh(x)\mathrm{d}g(x)\Biggl).
        \end{equation}
  \item[6).] Let $g, h$ be identities, then $\varphi\colon=f$. The isomorphic mean simplifies to the mean of the function. In this paper it is denoted by $\overline {{\displaystyle f}_{\scriptstyle {D}}}$, or $M_f$,
        \begin{equation}
            \overline {{\displaystyle f}_{\scriptstyle {D}}}=
            \frac{\int_{\scriptscriptstyle D}f(x)\mathrm{d}x}
            {\int_{\scriptscriptstyle D}\mathrm{d}x}~(=\isomeanvalue{f}{D}{\mathbb{V}x,\mathbb{V}y}).
        \end{equation}
      In the case $D$ is a close interval $[a,b]$,
        \begin{equation}
            \overline {{\displaystyle f}_{\scriptstyle {[a,b]}}}=\frac{1}{b-a}\int_a^bf(x)\mathrm{d}x.
        \end{equation}
  \item[7).] For monotone function $f$, its inverse function $f^{-1}$ is the dual-variable-isomorphic function of $f$ generated by mapping $f, f^{-1}$. Correspondingly the dual-variable-isomorphic mean of $f$ generated by $f, f^{-1}$ on a close interval is
      \begin{equation}
          \isomeanvalue{f}{[a,b]}{f,f^{-1}} = f\biggl(\frac{1}{f(b)-f(a)}\int_{f(a)}^{f(b)}f^{-1}(u)\mathrm{d}u\biggl).
      \end{equation}
\end{enumerate}
}
For convenience, while without confusions, we may use ``class I'' or ``Isomorphic mean class I'' for short name of Isomorphic mean class I of a function, and the similar for those of other classes later in this article, e.g. class II, class V.

For Theorem \ref{thm:IsoMeanInvar1}, especially when $g(x)=x$, then $M_f|_{m,h} = M_f|_{g,h} = M_f|_h$.

For Theorem \ref{thm:IsoMeanInvar2}, especially when $h(y)=y$, then $M_f|_{g,l} = M_f|_{g,h} = M_f|_g^{II}$.

For Corollary \ref{cor:IsoMeanInvar12}, especially when $g(x)=x$, $h(y)=y$, then $M_f|_{m,l} = M_f|_{g,h} = M_f$.
\subsection{Isomorphic mean class I of a function}\label{subsec:IsoMeanValFuncT1}
\begin{equation}
    \isomeanvalue{f}{[a,b]}{h}=M_f|_h=h^{-1}\biggl(\frac{1}{b-a}\int_a^bh\big(f(x)\big)\mathrm{d}x\biggl).
\end{equation}

As a special case, isomorphic mean class I also can be an equivalent derivation from isomorphic mean of numbers of $f(\xi_i)$, since in (\ref{equ:PartialSum}), let $g$ be identity mapping, then $\tau = f$. The class I can also be called the quasi-arithmetic mean of a function.

\subsubsection{Some properties of isomorphic mean class I}
\slthm{\label{thm:IsoMeanC1PrdInvariant} $\isomeanvalue{H_{ss}\big(f:k,C\big)}{[ka+C,kb+C]}{h}=\isomeanvalue{f}{[a,b]}{h}$.}
\slprf{
\begin{equation} \nonumber
\begin{split}
    \isomeanvalue{H_{ss}\big(f:k,C\big)}{[ka+C,kb+C]}{h} &=  h^{-1}\biggl(\frac{1}{(kb+C)-(ka+C)}\int_{ka+C}^{kb+C}h\big(f((u-C)/k)\big)\mathrm{d}u\biggl) \\
    &\overset{(x=(u-C)/k)}{=} h^{-1}\biggl(\frac{1}{b-a}\int_a^bh\big(f(x)\big)\mathrm{d}x\biggl).
\end{split}
\end{equation}
}
This means isomorphic mean class I is invariant with the function's horizontal scale and shift.

Isomorphic mean class I is not a homogeneous mean in general cases, as with constant $k$, generally
\begin{equation}
  h^{-1}\biggl(\frac{1}{b-a}\int_a^bh\big(kf(x)\big)\mathrm{d}x\biggl) \neq
      kh^{-1}\biggl(\frac{1}{b-a}\int_a^bh\big(f(x)\big)\mathrm{d}x\biggl),
\end{equation}
i.e. $M_{kf}|_h\neq kM_f|_h$.

There are some special cases or instances of it, as discussed below.

\subsubsection{Arithmetic mean of a function}
When $h\in \mathbb{V}y$,
\begin{equation}
    \isomeanvalue{f}{[a,b]}{h}=\overline {f(x)}
    =\frac{1}{b-a}\int_a^bf(x)\mathrm{d}x.
\end{equation}
It is the arithmetic mean of $f$ on $[a, b]$.

\subsubsection{Geometric mean of a function} \label{subsubsec:GeoMeanFunc}
Let $h(y)=\ln y$, $g$ be identity and $f(x)>0$ in (\ref{equ:PartialSum}), it turns into
\begin{equation} 
    M_{f(\xi_i)} =  \exp\biggl( \frac{1}{b-a}\sum_{i=1}^{n}{\ln f(\xi_i)\Delta x_i}\biggl) =  \sqrt[b-a]{\prod_{i=1}^{n}{f(\xi_i)^{\Delta x_i}}}.
\end{equation}
It is the weighted geometric mean of all $f(\xi_i)$. Especially the partitions being all equal it's the geometric mean of all $f(\xi_i)$,
\begin{equation} 
    M_{f(\xi_i)} = \sqrt[n]{\prod_{i=1}^{n}{f(\xi_i)}}.
\end{equation}
Then its limit form (\ref{equ:PartialSumLimit}) is deemed as positive function ``$f$'s geometric mean''. Therefore when $h(y)=\ln y, ~f(x)>0$, the following is called the geometric mean (value) of $f(x)$ on $(a, b)$
\begin{equation}
    \isomeanvalue{f(x)}{(a,b)}{\ln y}
    =e^{\frac{1}{b-a}\int_{a^{\scriptscriptstyle +}}^{b^{\scriptscriptstyle -}}\ln f(x)\mathrm{d}x}.
\end{equation}
And the following is called the geometric mean (value) of $f(x)$ on $[a, b]$
\begin{equation}
    \isomeanvalue{f(x)}{[a,b]}{\ln y}
    =e^{\frac{1}{b-a}\int_a^b\ln f(x)\mathrm{d}x}.
\end{equation}
It's obvious that geometric mean of a positive function is a homogeneous mean.

For a point $c\in[a,b]$, it holds that:
\begin{equation}
    {(\isomeanvalue{f}{[a,c]}{\ln y})}^{(c-a)} \cdot {(\isomeanvalue{f}{[c,b]}{\ln y})}^{(b-c)} = {(\isomeanvalue{f}{[a,b]}{\ln y})}^{(b-a)}. \label{eqn:GeoMeanFunMerge}
\end{equation}

Following are some special instances of geometric means on open intervals, in which the PVDM $\ln y$ is not defined on the $\inf M=0$($M$ is range of $f$). These only need an extra step of generalization by shrinking the domains $(a,b)$ of $f$ a little bit to $(a',b')$ so that $\inf M>0$, then taking the limit of the isomorphic mean with $a'\to a$ and/or $b'\to b$.
\begin{itemize}
  \item[(1).] The geometric mean of  $f(x)=x$ on $(0, b)$ : $\isomeanvalue{f(x)}{(0,b)}{\ln y}
        =e^{\frac{1}{b-0}\int_{0^{\scriptscriptstyle +}}^{b^{\scriptscriptstyle -}}\ln x\mathrm{d}x}=\frac be$.
  \item[(2).] The geometric mean of $f(x)=\sin x$ on $(0, \pi)$:
        $\isomeanvalue{\sin x}{(0,\pi)}{\ln y}
        =e^{\frac{1}{\pi}\int_{0^{\scriptscriptstyle +}}^{{\pi}^{\scriptscriptstyle -}}\ln \sin x\mathrm{d}x}.$
    Because the improper integral $\int_{0^{\scriptscriptstyle+}}^{{\pi}^{\scriptscriptstyle -}}\ln \sin x\mathrm{d}x=-\pi\ln2$ \cite{KOSMALA2}, such
    \begin{equation}
        \isomeanvalue{\sin x}{(0,\pi)}{\ln y}=e^{-\ln2}=\frac 12.
    \end{equation}
    It's easy to further conclude the geometric mean of $sin x$ on $(0, \frac{\pi}{2})$ is also $\frac 12$.
  \item[(3).] Function $tanx$ is unbounded on $(0, \frac{\pi}{2})$, however we consider the following as the generalized geometric mean for unbounded $tanx$ on $(0, \frac{\pi}{2})$, and we are able to work out the geometric mean value is 1.
    \begin{equation} 
       \isomeanvalue{\tan x}{(0,\pi/2)}{\ln y}
        =e^{\frac{2}{\pi}\int_{0^{\scriptscriptstyle +}}^{{\frac{\pi}{2}}^{\scriptscriptstyle -}}\ln \tan x\mathrm{d}x}=1.
    \end{equation}
  \item[(4).] A round is with diameter $d$, radius $r$. To compute the geometric mean of all parallel chords(e.g. in vertical direction): $\bar c$.
    The function of such chords can be written as $c=2 \sqrt{r^2-x^2}, (x\in(-r,r)$), such
    \begin{equation} 
        \bar c = \exp\bigl(\frac{1}{r-(-r)}\int_{(-r)^+}^{r^-}{\ln(2 \sqrt{r^2-x^2})\mathrm{d}x}\bigl) = \frac{4}{e}r = \frac{2}{e}d \approx 0.7358d.
    \end{equation}
\end{itemize}

\begin{figure}[htb]
\begin{center}
\includegraphics[keepaspectratio=true,scale=0.7]{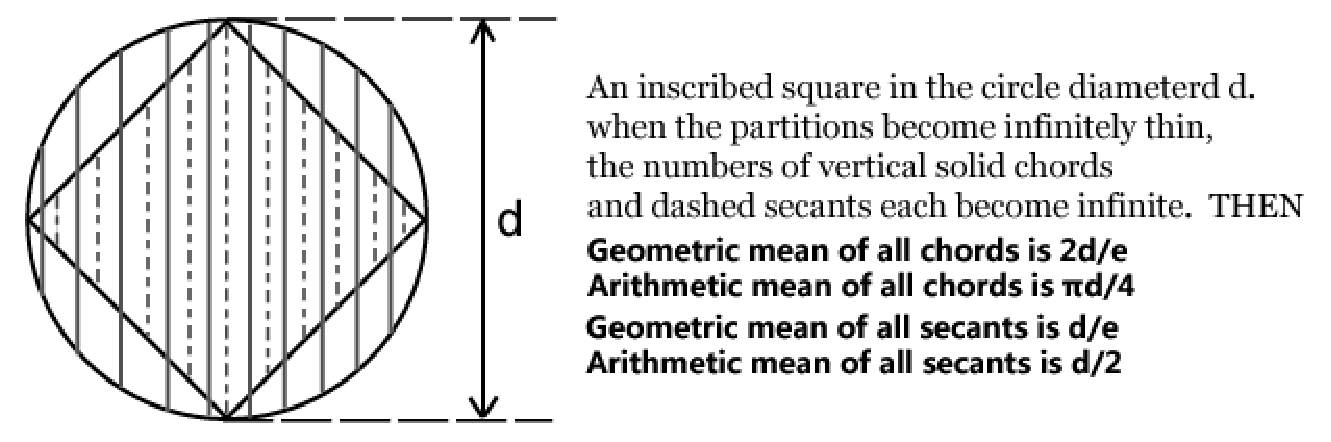}
\end{center}
\begin{illustration}{
\begin{center}
\label{illus:GeoArithRndSqr}Geometrical examples of geometric means of a function
\end{center}
}\end{illustration}
\end{figure}

Summarizing instance (1) and (4), and considering the homogeneity of geometric mean and its property by (\ref{eqn:GeoMeanFunMerge}) with instance (1), we have the following geometrical representation of some geometric means in Illustration \ref{illus:GeoArithRndSqr}.

\subsubsection{Harmonic mean of a function}
Let $h(y)=1/y$, $f(x)>0$ or $f(x)<0$,
\begin{equation}
    \isomeanvalue{f(x)}{[a,b]}{1/y}
    =\frac{b-a}{\int_a^b\frac{\mathrm{d}x}{f(x)}}
\end{equation}
is called the harmonic mean of $f(x)$ on $[a, b]$.

\subsubsection{Mean of power integral of a function}
Let $h(y)=y^p, ~(p\ne0)$, and $f(x)>0$ is continuous on $[a, b], ~b>a>0$,
\begin{equation}
    \isomeanvalue{f(x)}{[a,b]}{y^p}
    =\biggl(\frac{1}{b-a}\int_a^bf^p(x)\mathrm{d}x\biggl)^{\frac 1p},
\end{equation}
It is $f$'s $p$-order mean of power integral. Especially when $f(x)=x(p>-1,~ p\ne0, ~b>a>0)$,
\begin{equation}\label{equ:P-orderMean}
    \isomeanvalue{f(x)}{[a,b]}{y^p}
    =\biggl(\frac{b^{(p+1)}-a^{(p+1)}}{(p+1)(b-a)}\biggl)^{\frac 1p},
\end{equation}

\subsubsection{Function value on the midpoint of an interval}
When it happens that $f=h^{-1}$, then we know $f$ is monotone and $h=f^{-1}$, then
\begin{equation}
    \isomeanvalue{f(x)}{[a,b]}{f^{-1}}
    =f\biggl( \frac{1}{b-a}\int_a^bx\mathrm{d}x \biggl)=f\big(\frac{a+b}{2}\big).
\end{equation}
i.e. the monotone function's value on the midpoint of an interval is a special case of  isomorphic mean class I of the function.

\subsection{Isomorphic mean class II of a function}
 Thanks to the introduction of isomorphic frame and DVI function, isomorphic mean class II of a function of case 2) of Definition \ref{Def:IsoMeanOfFunSub} is the sibling of class I, though these two are quite different in their special forms.
\begin{equation}
    \isomeanvalueII{f}{[a,b]}{g}=M_f|_g^{II}=\frac{1}{g(b)-g(a)}\int_{g(a)}^{g(b)}f\big(g^{-1}(u)\big)\mathrm{d}u.
\end{equation}
If $g$ is differentiable,
\begin{equation}
    \isomeanvalueII{f}{[a,b]}{g}=M_f|_g^{II}=\frac{1}{g(b)-g(a)}\int_a^bf(x)\mathrm{d}g(x).
\end{equation}

Isomorphic mean class II is not invariant with H-scaleshift of the function in general cases, as opposed to that of class I, displayed in Theorem \ref{thm:IsoMeanC1PrdInvariant}.

\subsubsection{Some properties of isomorphic mean class II}
\paragraph{Homogeneity}
\slthm{$\isomeanvalueII{V_{ss}\big(f:k,C\big)}{[a,b]}{g}=k\isomeanvalueII{f}{[a,b]}{g}+C$.}
Proof omitted. This is again opposed to the property of class I.

\paragraph{Conjugation of 2 isomorphic means class II}
\slthm{ Let $g,~f$ be 2 strictly monotone, differentiable functions on interval $[a,b]$, and  $~A=f(a), B=f(b), C=g(a), D=g(b)$. There exists 4 isomorphic means class II: $~E=M_f|_g^{II}, ~F=M_g|_f^{II}$ ($~E\neq A, ~E\neq B, ~F\neq C,~F\neq\ D$),$~G=M_f|_f^{II}=(A+B)/2$, $~H=M_g|_g^{II}=(C+D)/2$. Then the following hold:
\begin{equation}
\begin{split}
  1).~&\big(\overrightarrow{AE} \colon \overrightarrow{EB}\big) \cdot \big( \overrightarrow{CF} \colon \overrightarrow{FD}\big) = 1 \label{thm:IsoMeanC2Conjug}\\
  2). ~&\overrightarrow{GE} \colon \overrightarrow{AB} = \overrightarrow{FH} \colon \overrightarrow{CD}.
\end{split}
\end{equation}
}

\begin{figure}[htb]
\begin{center}
\includegraphics[keepaspectratio=true,scale=1.1]{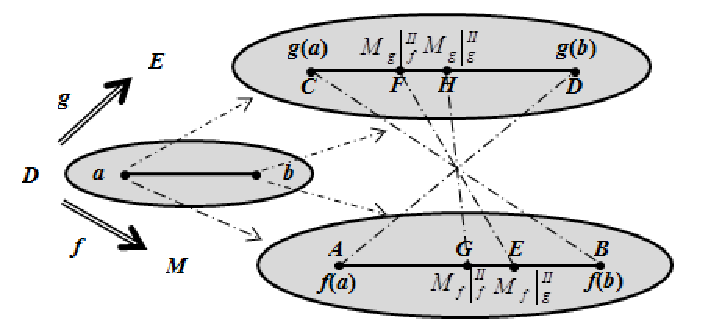}\\
$\big(\overrightarrow{AE} \colon \overrightarrow{EB}\big) \cdot \big( \overrightarrow{CF} \colon \overrightarrow{FD}\big) = 1$
\end{center}
\begin{illustration}{
\begin{center}
\label{illus:IsoMeanC2Conjug}Conjugation of 2 isomorphic means class II
\end{center}
}\end{illustration}
\end{figure}

\slprf{
\begin{equation} 
\begin{split}
        M_f|_g^{II} &=  \frac{1}{g(b)-g(a)}\int_a^bf(x)\mathrm{d}g(x)
         =  \frac{1}{g(b)-g(a)}\biggl( \big(f(x)g(x)\big)|_a^b - \int_a^bg(x)\mathrm{d}f(x)   \biggl) \\
        &=  \frac{1}{D-C}\big( B D -A C - F (B-A) \big)
        =  \frac{1}{D-C}\big( B (D-F) + A (F-C) \big)\\
        \Longrightarrow E &=  B \frac{\overrightarrow{FD}}{\overrightarrow{CD}} + A \frac{\overrightarrow{CF}}{\overrightarrow{CD}}
        \Longrightarrow E-A = B \frac{\overrightarrow{FD}}{\overrightarrow{CD}} + A \big( \frac{\overrightarrow{CF}}{\overrightarrow{CD}}-1 \big)
                            = B \frac{\overrightarrow{FD}}{\overrightarrow{CD}} - A \big( \frac{\overrightarrow{FD}}{\overrightarrow{CD}}\big).
        \nonumber
\end{split}
\end{equation}
It follows that $\overrightarrow{AE} \colon \overrightarrow{AB} = \overrightarrow{FD} \colon \overrightarrow{CD}$; Similarly $\overrightarrow{EB} \colon \overrightarrow{AB} = \overrightarrow{CF} \colon \overrightarrow{CD}$. Then it concludes with the result 1). The result 2) is a natural geometrical corollary of 1).}

Here the result 1) is said to be a kind of ``Conjugation of these 2 related isomorphic means class II($E$ and $F$)'' with exchanged original function and generator function. Illustration \ref{illus:IsoMeanC2Conjug} is the visual impression of the property when $f,g$ are of the same monotonicity.

\subsubsection{Simple special cases of isomorphic mean class II}\label{subsubsec:SpecCaseIsoMeanFuncT2}
\begin{enumerate}
  \item[(1).] Let $g\in \mathbb{V}x$,
    \begin{equation}
        \isomeanvalueII{f}{[a,b]}{\mathbb{V}x}=\overline {f(x)}
        =\frac{1}{b-a}\int_a^bf(x)\mathrm{d}x.
    \end{equation}
It is the arithmetic mean of $f(x)$ on $[a, b]$.
  \item[(2).] When $f(x)$ is monotone and continuous on $[a, b]$, and $g(x)=f(x)$,
    \begin{equation}
        \isomeanvalueII{f}{[a,b]}{g} = \isomeanvalueII{g}{[a,b]}{f}
        =\frac{1}{f(b)-f(a)}\int_{f(a)}^{f(b)}u\mathrm{d}u=\frac12\bigl(f(a)+f(b) \bigl),
    \end{equation}
    i.e. the arithmetic mean of the monotone function's values at the 2 endpoints of an interval is a special case of isomorphic mean class II. 
  \end{enumerate}

\subsubsection{A special case: With a classical mean value theorem}
\slthm{If $f\colon [a, ~b] \to \mathbb{R}$ is continuous and $g$ is an integrable function that does not change sign on $[a, b]$, then there exists $\xi$ in (a, b) such that
\begin{equation}
\begin{split}
  f(\xi)= \isomeanvalueII{f}{[a,b]}{G} ~~\biggl(=\int_a^b f(x) g(x)\mathrm{d}x /\int_a^bg(x)\mathrm{d}x\biggl) ~
         ~~(where~~G\in\int g).
\end{split}
\end{equation}
}
\slprf{With the conditions and the general case of well-known ``First mean value theorem for definite integrals'', there is a mean value $f(\xi)$ such
\begin{equation}
 \int_a^b f(x) g(x)\mathrm{d}x =f(\xi)\int_a^bg(x)\mathrm{d}x~~(\xi\in(a,b)).
\end{equation}
Let $G(x)=\int_a^x g(t)\mathrm{d}t+C$, then $G(x)$ is monotone and continuous, $G'(x)=g(x)$, thus
\begin{equation} 
\begin{split}
  f(\xi) =  \int_a^b f(x) g(x)\mathrm{d}x \div \int_a^b g(x)\mathrm{d}x = \frac{1}{G(b)-G(a)} \int_a^b f(x) \mathrm{d}G(x).
\end{split}
\end{equation}
Hence $f(\xi)$ is the isomorphic mean class II generated by the antiderivative of $g$.}
\slrem{While $f$ has not to be monotone, such $\xi$ has not to be unique, but the isomorphic mean $f(\xi)$ is unique.}
Especially when $g(x)=k,~G(x)=kx+C\in\mathbb{V}x$,
\begin{equation} 
\begin{split}
  f(\xi) = \isomeanvalueII{f}{[a,b]}{\mathbb{V}x} = \int_a^b kf(x)\mathrm{d}x \div \int_a^b k~\mathrm{d}x = \frac{1}{b-a} \int_a^b f(x) \mathrm{d}x
         =  \overline {{\displaystyle f}_{\scriptstyle {[a,b]}}}.
\end{split}
\end{equation}

\subsubsection{A special case: Elastic mean of a function}
  \sldef{The isomorphic mean class II of $f$ on $[a,b]~(b>a>0)$ generated by a logarithmic function e.g. $g(x)=\ln x(x>0)$,
    \begin{equation}\label{equ:logMeanFun}
        \isomeanvalueII{f}{[a,b]}{\mathbb{V}\ln x}
        =\frac{1}{\ln b-\ln a}\int_a^b\frac{f(x)}{x}~\mathrm{d}x,
    \end{equation}
is defined as the elastic mean of $f$ on $[a,b]$ in this paper.
}
\slrem{Any other logarithmic function is in $\mathbb{V}\ln x$.}
With $f(x)=x$,
\begin{equation} \label{equ:logMeanX}
        \isomeanvalueII{f}{[a,b]}{\ln x}
        =\frac{b-a}{\ln b-\ln a}.
\end{equation}

In Economics term, if $f$ is the elasticity \cite{TTElas} of another function $F$ such $f=xF'/F= (\mathrm{d}F/F)/(\mathrm{d}x/x)$, then it reflects the local relative change of $F$ against that of $x$, e.g. percentage change of demand (as an ``elastic'' response) against percentage change of price. Above $f(x)=x$ is the elasticity of $F(x)=ce^x (c > 0)$.

Let $K=F(b)\div F(a)$, and $k=b\div a$, it is easy to get
\begin{equation}
        \isomeanvalueII{f}{[a,b]}{\ln x}
        =\ln \big(F(b)\div F(a)\big) \div \ln(b \div a)=\log_{~k} K. 
\end{equation}
On the other hand, let $M=F(x+\mathrm{d}x)\div F(x)$, $m=(x+\mathrm{d}x)\div x$, then
\begin{equation}
\begin{split}
        f(x) &=(\mathrm{d}F/F)/(\mathrm{d}x/x)=\mathrm{d}\ln F/\mathrm{d}\ln x \approx \Delta \ln F/\Delta \ln x \\
            &=(\ln F(x+\mathrm{d}x)-\ln F(x))\div(\ln(x+\mathrm{d}x)-\ln x)=\log_m M,
\end{split}
\end{equation}
i.e. $f$ is a ``micro-logarithm'' of 2 ``micro-multiplications'' for every local $x$, while $\isomeanvalueII{f}{[a,b]}{\ln x}$ is the logarithm of 2 overall multiplications over $[a,b]$.

Therefore the elastic mean of $f$ is actually the ``average of the elasticity $f$''. It is very similar to arithmetic mean of a function $f$ computed via the Newton-Leibniz formula:
\begin{equation}
        M_f = {\big(F(b)-F(a)\big)}\div {(b-a)},
\end{equation}
where $F$ is now the anti-derivative of $f$. 

\subsubsection{Instances of isomorphic mean class II}
\paragraph{Elastic mean of \texorpdfstring{$\tan x$\\}{tan x}}
A special instance of elastic mean is about unbounded $tanx$ on $(0, \frac{\pi}{2})$, as following limiting form:
\begin{equation}
\begin{split}
   \isomeanvalueII{\tan x}{(0,\pi/2)}{\ln x} &= \lim_{x\to 0,y\to\frac{\pi}{2}} \frac{1}{\ln \frac{\pi}{2}-\ln x}\int_0^y \frac{\tan(t)}{t}\mathrm{d}t.
\end{split}
\end{equation}
The numerator part is an improper definite integral, and the denominator part is approaching $+\infty$. We transform it with $~x=\rho\cos\theta,~y=\frac{\pi}{2}+\rho\sin\theta ~~(\rho>0, ~-\frac{\pi}{2}<\theta<0)$, and apply L~'Hopital's ~rule for twice in the following:
\begin{equation}
\begin{split}
   \isomeanvalueII{\tan x}{(0,\pi/2)}{\ln x} &= \lim_{\rho \to 0} \frac{1}{\ln \frac{\pi}{2}-\ln (\rho\cos\theta)}\int_0^{\frac{\pi}{2}+\rho\sin\theta} \frac{\tan(t)}{t}\mathrm{d}t \\
    &= \lim_{\rho \to 0}\frac{1}{-\cos\theta\frac{1}{\rho\cos\theta}}\cdot\frac{\tan(\frac{\pi}{2}+\rho\sin\theta)}{\frac{\pi}{2}+\rho\sin\theta}\cdot \sin\theta
    = \frac{2}{\pi}.
\end{split}
\end{equation}

\paragraph{Elastic mean of power function\texorpdfstring{\\}{}}
Let $f(x)=x^p, p\ne0, b>a>0, x\in [a, b], g(x)=\ln x(x>0)$,
\begin{equation} 
\begin{split}
    \isomeanvalueII{f}{[a,b]}{\ln x} =  \frac{1}{\ln b-\ln a}\int_{\ln a}^{\ln b}f(e^u)\mathrm{d}u 
    =  \frac{f(b)-f(a)}{\ln f(b)-\ln f(a)},
\end{split}
\end{equation}
which is the logarithmic mean of $f(a)$ and $f(b)$.

\paragraph{Power function's isomorphic mean class II generated by \texorpdfstring{$1/x(x>0)$\\}{1/x(x>0)}}
Let $f(x)=x^p, p\ne1, b>a>0, x\in [a, b], g(x)=1/x(x>0)$,
\begin{equation} 
\begin{split}
    \isomeanvalueII{f}{[a,b]}{1/x} = \frac{1}{(1/b)-(1/a)}\int_{1/a}^{1/b}(1/u)^p\mathrm{d}u 
    = \frac{ab(b^{p-1}-a^{p-1})}{(p-1)(b-a)}.
\end{split}
\end{equation}
When $p=2,~\isomeanvalueII{f}{[a,b]}{1/x}=ab$; when $p=3,~\isomeanvalueII{f}{[a,b]}{1/x}=\frac12ab(a+b)$.

While in the case $p=1, ~\isomeanvalueII{f}{[a,b]}{1/x}=\frac{1}{(1/b)-(1/a)}\int_{1/a}^{1/b}(1/u)\mathrm{d}u=\frac{ab(\ln b-\ln a)}{(b-a)}$, It is the product of $a, b$ divided by the logarithmic mean of $a, b$, it is also the limit of above $\frac{ab(b^{p-1}-a^{p-1})}{(p-1)(b-a)}$ when $p\to 1$.

\subsection{Isomorphic mean class III  \& IV of a function }
Isomorphic mean class III \& class IV are just the general forms of the Definition, being the combined form of class I \& class II. Their properties are mainly covered by previous sections.

\subsubsection{A special case of class III}
\slthm{\label{thm:IsoMeanC3ofIden} An isomorphic mean class III of $f(x)=x$ on [a,b] generated by $g$ equals the isomorphic mean of $a,~b$ generated by $g$ (the generalized $g$-mean).
}
\slprf{
\begin{equation}
    \isomeanvalueIII{x}{[a,b]}{g} =g^{-1}\biggl(\frac{1}{g(b)-g(a)}\int_{g(a)}^{g(b)}g\bigl(g^{-1}(u)\bigl)\mathrm{d}u\biggl)
    =g^{-1}\big(\frac{g(a)+g(b)}{2}\big).
\end{equation}}

\subsection{Isomorphic mean class V} 
Isomorphic mean class V in (\ref{equ:IsoMeanT5}) can be deemed as a special mean of a single variable, e.g. in the form of
\begin{equation} \nonumber
    \isomeanvalue{x}{[a,b]}{g,h}=h^{-1}\Biggl(\frac{1}{g(b)-g(a)}\int_a^bh(x)\mathrm{d}g(x)\Biggl).
\end{equation}
While traditionally without isomorphic frame in mind, it is not so worthwhile to discuss, since on $[a,b]$ $y=x$ always has a mean value $\frac{1}{2}(a+b)$. (Recalling Section \ref{subsubsec:natureCompIsoMean}, in article \cite{LIUY} the so-called extended convexity of $y=x$ is somehow meaningful too.) The class V is generally not a mean value of an ordinary function. It is further sub-classifiable with $g\in \mathbb{V}x$ or $h\in \mathbb{V}y$, corresponding to class I or class II of $f(x)=x$, i.e. of a single variable, which are the closest concepts to so-called ``class 0'': the isomorphic mean of numbers.

\subsubsection{Composite class V}
\sldef{With a strictly monotone and continuous $f$, and applicable $g,h$,
\begin{equation}
M_x|_{g,H}=f^{-1}(M_f|_{g,h})
\end{equation}
is called a composite isomorphic mean class V generated by $g,h,f$, where $H:=h\circ f$.}
Above $H$ is strictly monotone and continuous. Composite class V is a special case of class V, but from the view of a normal class V of $t=M_x|_{g,H}$, the PVDM $H$ can be decomposed to $h\circ f$ whereby class V can be related to mean values of more monotone functions, i.e.
$f(t)=M_f|_{g,h}.$

Both class V and composite class V are useful when relating to other types of means in mathematics. See Section \ref{sec:RelationIsoMeanCauchy}.

\subsubsection{Generation of bivariate means by class V}\label{sec:GenBivarMean}
For a close interval $[a,b]$ and applicable $g,h$(and a monotone $f$ as with composite class V), the $M_x|_{g,h}$ (or $M_x|_{g,H}$) is clearly a sort of mean value of $a,b$. Moreover with the ``property of symmetry with endpoints of interval'', such bivariate mean is symmetric. Below are 2 examples among possible others.

\paragraph{Bivariate means regarding trigonometric functions\texorpdfstring{\\}{}}
By choosing $g(x)=\sin x, ~h(y)=\cos y$, $[a,b]\subseteq[0,\pi/2]$,
\begin{equation}
\begin{split}
    \isomeanvalue{x}{[a,b]}{\sin x, \cos y} &= \arccos \biggl( \frac{1}{\sin b- \sin a} \int_a^b \cos x (\sin x)'\mathrm{d}x  \biggl) \\
        &= \arccos \biggl( \frac{b-a+\sin b\cos b-\sin a\cos a}{2(\sin b- \sin a)} \biggl).
\end{split}
\end{equation}
While $g,h$ are exchanged,
\begin{equation}
\begin{split}
    \isomeanvalue{x}{[a,b]}{\cos x, \sin y} &= \arcsin \biggl( \frac{1}{\cos b- \cos a} \int_a^b \sin x (\cos x)'\mathrm{d}x  \biggl) \\
        &= \arcsin \biggl( \frac{a-b+\sin b\cos b-\sin a\cos a}{2(\cos b- \cos a)} \biggl).
\end{split}
\end{equation}

\paragraph{A class of quasi-Stolarsky means\texorpdfstring{\\}{}}
Another example of bivariate mean class generated by the class V, in the case $g(x)=x^p(x>0, ~p\ne0), ~h(y)=y^q(y>0, q\ne0), a>0,b>0,a\ne b$. It's denoted and formulated by:
\begin{equation} \label{equ:quasiStolarsky}
    Q_{p,q}(a,b)= \biggl( \frac{p(b^{p+q}-a^{p+q})}{(p+q)(b^p-a^p)}\biggl)^{1/q}.
\end{equation}

It is very similar to the derivation of the Stolarsky means from the ``Cauchy's extended mean value theorem''(\cite{LEASCHOExtM2}, pp207) by a pair of power functions. This class has also different special cases(including limits of $Q_{p,q}(a,b)$ when $p\to0$ and/or $q\to0$, which are actually equivalent results as with replaced $g(x)=\ln x$ and/or $h(y)=\ln y$):
\begin{equation} \label{equ:quasiStolarsky2}
    \begin{split}
        Q_{p,q}(a,b)=\left\{
        \begin{aligned}
        &\big(\frac{a^p+b^p}{2}\big)^{1/p} &(p=q,pq\ne0), \\
        &\sqrt{ab} &(p=q=0), \\
        &\biggl(\frac{b^p-a^p}{p(\ln b-\ln a)}\biggl)^{1/p} &(p+q=0,pq\ne0), \\
        &\biggl(\frac{b^q-a^q}{q(\ln b-\ln a)}\biggl)^{1/q} &(p=0,q\ne0), \\
        &\exp\biggl(\frac{b^p\ln b-a^p\ln a}{b^p-a^p}-\frac{1}{p}\biggl) &(q=0,p\ne0), \\
        &\frac{2(a^2+ab+b^2)}{3(a+b)} &(p=2,q=1), \\
        &\sqrt[3]{a\cdot\frac{a+b}2\cdot b} &(p=-1,q=3). \\
        \end{aligned} \right.\\
    \end{split}  
   \end{equation}

Details of derivations and proofs are omitted. It features the class of power mean as its well-balanced children forms. Though conversion of this form to Stolarsky mean is easy by substituting power $s=p+q$, but from the perspective of isomorphic mean, it's the very balanced form with respect to $p$ and $q$.

\subsection{Isomorphic mean class VI,VII}
The case 6) of Definition \ref{Def:IsoMeanOfFunSub} is just the mean of $f$, as ``class VI'' of isomorphic mean. It is not further discussed here. For instance of case 7) as ``class VII'', if $f(x)=x^a (x\in [0, c], a>0)$ then $  \isomeanvalue{f}{[0,c]}{f,f^{-1}} = (\frac{a}{a+1})^ac^a$, which coefficient $(\frac{a}{a+1})^a$ is approaching $\frac{1}{e}$ when $a$ is approaching $+\infty$.

\section{\label{sec:RelationIsoMeanCauchy}Isomorphic mean of a function vs Cauchy mean value}

\subsection{About Cauchy mean value}
\slthm{\label{thm:CauchyMeanVal}Cauchy's mean-value theorem states that (\cite{LOSONLCAUCHYCOMP}, pp12): If $f$, $g$ are continuous real functions on $[x_1,x_2]$ which are differentiable in $(x_1,x_2)$, and $g'(u)\ne 0$ for $u\in(x_1,x_2)$, then there is a point  $t\in(x_1,x_2)$ such that
\begin{equation}
\frac{f'(t)}{g'(t)}=\frac{f(x_2)-f(x_1)}{g(x_2)-g(x_1)}.
\end{equation}}

Note such $t$ does not have to be unique in general cases(e.g. an definition in \cite{HILLE} says there is at least one such $t\in(x_1,x_2)$ are satisfactory). Hence none of $f(t)$,$g(t)$,$f'(t)$,$g'(t)$ has to be unique, i.e. there is not uniquely defined mean value of function(s) either.

For such $t$ to be unique, in \cite{LOSONLCAUCHYCOMP} there is further restriction, whereby there comes the definition of ``Cauchy mean value of two numbers'':
\sldef{\label{Def:CauchyMeanof2Number}Assuming now (with Theorem \ref{thm:CauchyMeanVal}) that $f'/g'$ is invertible we get
\begin{equation}
t=\big(\frac{f'}{g'}\big)^{-1}\bigg(\frac{f(x_2)-f(x_1)}{g(x_2)-g(x_1)}\bigg).
\end{equation}
This number is called the Cauchy mean value of the numbers $x_l, x_2$ and will be denoted by $t= D_{fg}(x_1, x_2)$.}
In this paper, we further say the mean value of such 2 numbers ``is generated by $f,g$''. (And in \cite{LOSONLCAUCHYCOMP}, further covered is a generalized form of above to Cauchy mean of $n$ numbers: $D_{fg}(x_1, x_2,...,x_n)$). But even with Definition \ref{Def:CauchyMeanof2Number}, neither $f(t)$ nor $g(t)$ can be deemed as well-defined mean value of a function, as they are symmetrically, temporarily depending on each other.

The Cauchy mean value class and certain classes of isomorphic means can be converted to each other, with some criteria. A reasonable point is that since Cauchy mean value deals with 2 functions while isomorphic mean does with 3, the corresponding isomorphic mean is related to or of class V, which is with simplest $f(x)=x$.

\subsection{Conversion from Cauchy mean value to class V}
To get a class V from an ordinary Cauchy mean value, below supporting theorem is needed:
\slthm{(Extended Darboux's Theorem:) If on $[a,b]$ $f,g$ are both differentiable and $g'(x)\ne0$, then on $(a,b)$ $f'(x)/g'(x)$ can take any value between $f'(a)/g'(a)$ and $f'(b)/g'(b)$.
}

Then we have the following theorem:
\slthm{\label{thm:CauchyMean2IsoMean}Let $f$, $g$ be differentiable on $[a,b]$ and $g'$ is Riemann integrable, $\forall x\in(a,b)~g'(x)\ne 0$ and $g(b)\ne g(a)$. If $h:=f'/g'$ is invertible on $[a,b]$, then the Cauchy mean value of $a,b$ generated by $f,g$ is the isomorphic mean class V on $[a,b]$ generated by $g,h$. i.e.
The Cauchy mean value $t\in(a,b)$ is such that
\begin{equation}
\begin{split}
    t&=\big(\frac{f'}{g'}\big)^{-1}\bigg(\frac{f(b)-f(a)}{g(b)-g(a)}\bigg) = D_{fg}(a,b)\\
     &=~ h^{-1}\Biggl(\frac{1}{g(b)-g(a)}\int_a^bh(x)g'(x)\mathrm{d}x\Biggl) = \isomeanvalue{x}{[a,~b]}{g,h}.
\end{split}
\end{equation}
}
\slprf{Referring to Definition \ref{defin:IsoMeanOfFun} and \ref{Def:IsoMeanOfFunSub}, we claim that both $h,g$ are continuous and strictly monotone on $[a,b]$, and $y=x$ has the isomorphic mean. \\
\indent (i) According to the Extended Darboux's Theorem, $h=f'/g'$ has the intermediate value property(IVP): $\forall c, \forall d\in [a,b]~(c<d)$ and every $w$ between $h(c)$ and $h(d)$~(here $h(c)\ne h(d)$ due to invertible $h$), $\exists \xi\in(c,d)$ such that $h(\xi)=w$. Moreover because it's invertible, $h$ is injective. With these 2 properties, It's easy to prove that $h$ is further strictly monotone, and again with IVP $h$ is continuous on $[a,b]$.\\
\indent (ii) $g$ is continuous due to $g$ is differentiable. Since $g'$ is non-zero, according to Darboux's Theorem $g'$ is either all positive or all negative on $(a,b)$, and $g'_+(a)$ or $g'_-(b)$ also has the same sign as they are non-zero for $f'_+(a)/g'_+(a)$ or $f'_-(b)/g'_-(b)$ to exist. Such $g$ is strictly monotone.\\
\indent (iii) As $y=x$ is continuous, according to Theorem \ref{thm:IsomeanExist1}, $\isomeanvalue{x}{[a,b]}{g,h}$ exists.\\
\indent (iv) The Cauchy mean value $t$ is (uniquely) existing as $h$ is invertible, and these 2 mean values are easily seen equal when $g'$ is Riemann integrable(such $hg'$ also Riemann integrable), since according to the fundamental theorem of calculus:
\begin{equation}
  \int_a^bh(x)g'(x)\mathrm{d}x=\int_a^bf'(x)\mathrm{d}x=f(b)-f(a).
\end{equation}
\indent This ends the proof of the theorem.}
\slrem{With Theorem \ref{thm:CauchyMean2IsoMean}, practically if $h:=f'/g'$ is continuous and strictly monotone, then it's invertible, and $\isomeanvalue{x}{[a,~b]}{g,h}$ exists.}

For example: With Theorem \ref{thm:CauchyMean2IsoMean}, let 2 generator functions be $f(x)=\ln x$, $g(x)=x$, $b>a>0$,
\begin{equation}
D_{fg}(a,b)=\frac1{\frac{f(b)-f(a)}{g(b)-g(a)}}=\frac{b-a}{\ln b-\ln a}.                      
\end{equation}
In order to get class V, $h(x)=f'(x)/g'(x)=1/x$. To verify:
\begin{equation}
\isomeanvalue{x}{[a,b]}{g,h}=\frac1{\frac{1}{b-a}\int_a^b\frac{1}{x}\mathrm{d}x}=\frac{b-a}{\ln b-\ln a}=D_{fg}(a,b).
\end{equation}

\subsection{Conversion from class V to Cauchy mean value}
\slthm{\label{thm:IsomeanV2CauchyMean} Let $f$ be a strictly monotone, continuous and bounded function on $[a,b]$ with range $M$; Let $g$ be strictly monotone and differentiable on $[a,b]$, $g'$ is Riemann integrable and $\forall x\in(a,b)~g'(x)\ne 0$; Let $h$ be strictly monotone continuous on $M$ and $\gamma:=h\circ f$. Then there are following equivalences:\\
\indent1). 
\begin{equation}\label{equ:ConvIsoC5toCauchyMean}
\begin{split}
    t&=\isomeanvalue{x}{[a,b]}{g,\gamma}=f^{-1}\Biggl(h^{-1}\biggl(\frac{1}{g(b)-g(a)}\int_a^bh(f(x))g'(x)\mathrm{d}x\biggl)\Biggl)\\
     &=~ \gamma^{-1}\biggl(\frac{l(b)-l(a)}{g(b)-g(a)}\biggl)
     =\big(\frac{l'}{g'}\big)^{-1}\biggl(\frac{l(b)-l(a)}{g(b)-g(a)}\biggl)=D_{lg}(a,b).\\
     \big(~&where~~ l\in\int(\gamma~g')=\int((h\circ f) g').~\big)
\end{split}
\end{equation}
\indent2). 
\begin{equation}
    \isomeanvalue{f(x)}{[a,~b]}{g,h}=f(t)=f(D_{lg}(a,b)).
\end{equation}
Especially when $f$ is identity,
\begin{equation}
    \isomeanvalue{x}{[a,b]}{g,h}=t=D_{lg}(a,b), ~~where~~ l\in\int(h g').
\end{equation}
}
\slprf{With such defined $g,h,f$ and derived ~$l$, the $f^{-1}(\isomeanvalue{f(x)}{[a,b]}{g,h})$ and $D_{lg}(a,b)$ exist and they are equal as proved by the inline formulae (\ref{equ:ConvIsoC5toCauchyMean}), and further 2) is true.}

\slrem{Let $l(x):=\int_a^x (h(f(x))\cdot g'(x)dx+C, ~l_2(u):=\int_{g(a)}^u h(f(g^{-1}(u)))du+C$, then $l_2(u)=l\circ (g^{-1})(u)$. According to the Lagrange's differential mean value theorem, there is a $\xi\in(\min\{g(a),g(b)\},\max\{g(a),g(b)\})$ such that:
\begin{equation}
\begin{split}
    &l'_2(\xi)=\frac{l_2(g(b))-l_2(g(a))}{g(b)-g(a)} \\
     \overset{t=g^{-1}(\xi)}{\Longrightarrow} &\frac{l'(t)}{g'(t)}=\frac{l(b)-l(a)}{g(b)-g(a)}.
\end{split}
\end{equation}
This means: The composite isomorphic mean class V (by $g,h,f$) $t$ is the image of the Lagrange mean value $\xi$ of $g(a),~g(b)$(``generated by $\int\varphi(f:g,h)$'') under the inverted bijection of $g$.
}
Conversion examples:
\begin{enumerate}
  \item[(1).] Conversion of a geometric mean into a Cauchy mean value.\\
Geometric mean of a positive monotone $y=f(x)$, is a special case of class I where $g(x)=x$, $h(y)=\ln y$. With Theorem \ref{thm:IsomeanV2CauchyMean}, $l\in\int h(f(x))g'(x)\mathrm{d}x=\int \ln f(x)\mathrm{d}x$. Then
    \begin{equation}
    D_{lg}(a,b)=\bigl(\frac{l'}{g'}\bigl)^{-1}\bigl(\frac{l(b)-l(a)}{g(b)-g(a)}\bigl)=f^{-1}\Biggl(\exp\biggl(\frac{1}{b-a}\int_a^b\ln f(x)\mathrm{d}x\biggl)\Biggl).
    \end{equation}
  \item[(2).] Conversion of an elastic mean into a Cauchy mean value.\\
Elastic mean of a monotone $y=f(x)$, is a special case of class II where $g(x)=\ln x$, $h(y)=y$. With Theorem \ref{thm:IsomeanV2CauchyMean}, $l\in\int h(f(x))g'(x)\mathrm{d}x=\int (f(x)/x)\mathrm{d}x$.
    \begin{equation}
    D_{lg}(a,b)=\bigl(\frac{l'}{g'}\bigl)^{-1}\bigl(\frac{l(b)-l(a)}{g(b)-g(a)}\bigl)=f^{-1}\biggl(\frac{1}{\ln b-\ln a}\int_a^b\frac{f(x)}{x}\mathrm{d}x\biggl).
    \end{equation}
\end{enumerate}

\subsection{\label{sec:DifferenceIsoCauchy}Differentiations between isomorphic mean and Cauchy's mean value theorem}

\subsubsection{About their coverage and categorical intersection}
Since the IVP of isomorphic means applies to ordinary functions, while IVP of the Cauchy's mean value theorem applies to derivable functions, and in view of Theorem \ref{thm:CauchyMean2IsoMean} and Theorem \ref{thm:IsomeanV2CauchyMean}, we claim that (i). Whenever there is a Cauchy mean value, there is a class V(only of derivable DMs). (ii). Whenever there is an isomorphic mean, it's not ensured there is a convertible Cauchy mean value by our ready Theorem.

Following are some examples of isomorphic means that can not convert to Cauchy mean values.
\begin{itemize}\setlength{\itemsep}{-0.2em}
  \item[(1).] In Section \ref{subsec:ExiCompatIsoMean} case (1), the function $f$ taking value of either 1 or 3, has jump discontinuity, though the isomorphic mean is 2. With reference to Theorem \ref{thm:IsomeanV2CauchyMean}, we attempt an $l\in\int((h\circ f)\cdot g')=\int 2(f(x)+2)\mathrm{d}x$, $a=0,b=2$. It follows that $(l(b)-l(a))/(g(b)-g(a))=4$, but $l'(x)/g'(x)=h(f(x))=f(x)+2$ which value is either 3 or 5 on $[a,b]$, never being 4. This is simply because $f$ being not continuous at $x=1$ makes $l$ being not differentiable everywhere on $(a,b)$, therefore some criteria for Theorem \ref{thm:CauchyMeanVal} are not met. Such here the Cauchy mean value of 0,~2 even does not exist.
  \item[(2).] $f(x)=\sin x (x\in[0,\pi])$ being not monotone, has a mean value $2/\pi$, which is a special case of isomorphic mean (class VI). $(l(b)-l(a))/(g(b)-g(a))=2/\pi$, but $l'(x)/g'(x)=h(f(x))=f(x)$ is not invertible, therefore the Cauchy mean value does not exist (uniquely). However in this case Theorem \ref{thm:CauchyMeanVal} holds.
  \item[(3).] $f(x)\equiv A$, such with any applicable $g,h$ its isomorphic mean is $A$. It's easy to check Theorem \ref{thm:CauchyMeanVal} holds for every $x$ with $l,g$, but the Cauchy mean value does not exist either.
\end{itemize}

The IVP of an isomorphic mean primarily gets an unique intermediate value in-between of the function values, while the IVP of Cauchy's mean value theorem primarily non-unique intermediate values within the function's domain. The latter calls for an unique value in the domain by further restriction of $f'/g'$ being invertible, whereas the former calls for an unique value within the domain by specially letting $f(x)=x$, then there is an intersection of the two categories, i.e. convertible instances between the isomorphic means class V and the Cauchy mean values.

\textbf{In summary, isomorphic means of a function cover broader range of unique mean values of an ORDINARY function.} Only the class V can be matched by the Cauchy mean value as a special derivation of the Cauchy's mean value theorem in term of uniqueness, whereas more other classes of isomorphic means conform to the ubiquitous theorem. 

\subsubsection{About generator functions and identifications of unique means}
For an isomorphic mean, $f$'s bonding on $\{g,h\}$ makes the DMs $g,h$ not in the equal or symmetric statuses with the focused $f$, such the isomorphic mean is prominently an extended unique mean value of $f$ as of good identification. 

On the other hand for a Cauchy mean value there are only 2 generator functions $f,g$. Though $f(t)$ is somehow a kind of mean generated by $g$ and vise versa for $g(t)$,$f$, rather $f,g$ are symmetrical in statuses meanwhile $f'/g'$ being monotone is a requirement of mixed nature. Therefore the Cauchy mean value $t$ and its corresponding generator functions' values $f(t)$ or $g(t)$ are not well identifiable as being mean values of a specific function.

Cauchy mean value lacks of information of 1 ordinary function to match the complexity of a class IV, but only to its best it matches a special class V.

\textbf{In summary, isomorphic means have better identifications as being well-defined unique extended means of a specific function.} It is with such good identifications, that isomorphic means also have had diversified classifications possibly.

\subsubsection{Conclusion}
Generally speaking: Isomorphic mean seems a concept of better origin and perspective, better identification and classification, more coverage and more natural generalization.

However certain conversions between these 2 types are useful. For examples,
\begin{itemize}\setlength{\itemsep}{-0.2em}
  \item[(1).]With the quasi-Stolarsky formula (\ref{equ:quasiStolarsky2}), it's easy to know the power mean class is a Cauchy mean value generated by $f(x)=x^{2p}, ~g(x)=x^p$.
  \item[(2).]The solutions to comparisons of certain isomorphic means can be easily based on those to the comparison of Cauchy mean values, which applications are whereby made richer.
\end{itemize}

\section{The comparison problems of \texorpdfstring{\\}{}isomorphic means of a function}
We will discuss 2 types of the comparison methods of isomorphic means of a function: the methods derived from comparison methods of Cauchy mean values, and the methods by help of monotonicity and convexity conditions.

\subsection{Comparison methods of class V scenario derived from \texorpdfstring{\\}{}comparison of Cauchy mean values}
The Theorem 1 and Theorem 2 of Losonczi \cite{LOSONLCAUCHYCOMP} are the necessary conditions and sufficient conditions respectively for the functions $f,g,F,G$ such the comparison inequality(\cite{LOSONLCAUCHYCOMP} (2))
\begin{equation} \label{equ:CauchyMeanCompLos}
D_{fg}(x_1,x_2,...,x_n)\leq D_{FG}(x_1,x_2,...,x_n)
\end{equation}
holds, where $n\geq2$ is fixed.

The prerequisites for both theorems are:
$I$ is real interval, $\varepsilon_n(I)$ is the set of all pairs $(f,g)$ of functions $f,g:I\to\mathbb{R}$ satisfying
the following conditions:\\
(i) $f,g$ are $n$-times differentiable on $I$,\\
(ii) $g^{(n-1)}(u)\ne0$ for $u\in I$,\\
(iii) the (first) derivative of $(f^{(n-1)}/g^{(n-1)})$ is not zero on $I$.\\
And there are notations: $\bar f=f^{(n-1)},\bar g=g^{(n-1)},\bar F=F^{(n-1)},\bar G=G^{(n-1)},~h=\bar f/\bar g,~H=\bar F/\bar G $.

Losonczi's necessary conditions for (\ref{equ:CauchyMeanCompLos}) is quoted:
\slthm{\label{thm:Los1} (\cite{LOSONLCAUCHYCOMP}, pp15, Theorem 1)
Suppose that the functions $f,g,F,G:I\to\mathbb{R}$ are $n+1$ times continuously differentiable and $(f,g),(F,G)\in \varepsilon_n(I)$. Then the inequality
\begin{equation}
\frac{h''(x)}{h'(x)}+2\frac{\bar{g}~'(x)}{\bar{g}(x)} \leq \frac{H''(x)}{H'(x)}+2\frac{\bar{G}~'(x)}{\bar{G}(x)}~~ (x\in I)
\end{equation}
is necessary for (\ref{equ:CauchyMeanCompLos}) to hold.
}
According to remark 1 of \cite{LOSONLCAUCHYCOMP}, these necessary conditions' proof assumes $n$ values are such that $x_2=...=x_n$ are near $x_1$. This is especially meaningful for $n=2$ cases.

Losonczi's sufficient conditions for (\ref{equ:CauchyMeanCompLos}) is read as:
\slthm{\label{thm:Los2} (\cite{LOSONLCAUCHYCOMP}, pp18, Theorem 2)
Suppose that $(f,g),(F,G)\in \varepsilon_n(I)$. Then the inequality
\begin{equation}
\frac{h(u)-h(v)}{h'(v)}~\frac{\bar{g}(u)}{\bar{g}(v)} \leq \frac{H(u)-H(v)}{H'(v)}~\frac{\bar{G}(u)}{\bar{G}(v)}~~ (u,v\in I)
\end{equation}
is sufficient for (\ref{equ:CauchyMeanCompLos}) to hold.
}

To make them use for composite class V scenario, we consider the $n=2$ basic cases of both theorems. Then by coincidences or by internal relations, the above important parametric function $h:=\bar f / \bar g$ is right our PVDM $h:=f'/g'$ with Theorem \ref{thm:CauchyMean2IsoMean} or $\gamma:=l'/g'$ with Theorem \ref{thm:IsomeanV2CauchyMean}. This makes these 2 theorems literally much orientated to and ready for isomorphic means.

Our comparison problems are: To find necessary conditions and sufficient conditions for functions $f,g,h,G,H$ such the inequality
\begin{equation}\label{equ:IsomeanC5Comp}
\isomeanvalue{f(x)}{[a,b]}{g,h}\leq \isomeanvalue{f(x)}{[a,b]}{G,H}
\end{equation}
holds, where $f$ is continuous and invertible on $[a,b]$($a\ne b$). For valid comparison problems, all isomorphic means of a function in various context of this section are assumed existing with an applicable isomorphic frame.

We have following derived theorems:
\subsubsection{Necessary conditions}
\slthm{\label{thm:IsomeanC5CompNess}Let
(i) $f:[a,b]\to M$ be a strictly monotone, differentiable and bounded function;
(ii) $g,G$ be 2 times continuously differentiable bijections on $[a,b]$, on which $~g'\ne 0,~G'\ne 0$;
(iii) $h,H$ be 2 times continuously differentiable bijections on $M$, on which $~h'\ne 0,~H'\ne 0$. Then
\begin{equation} \label{equ:IsomeanVNecess}
\begin{split}
\frac{\gamma''(x)}{\gamma'(x)}+2\frac{g''(x)}{g'(x)} \leq(\geq) \frac{\Gamma''(x)}{\Gamma'(x)}+2\frac{G''(x)}{G'(x)}\\
\big(x\in[a,b],~~f~is~increasing(decreasing)\big)
\end{split}
\end{equation}
is necessary for $\isomeanvalue{f(x)}{[a,b]}{g,h}\leq \isomeanvalue{f(x)}{[a,b]}{G,H} $(\ref{equ:IsomeanC5Comp}) to hold, where $\gamma:=h\circ f$, $\Gamma:=H\circ f$.
}
\slprf{Such required $f,g,h,G,H$ make $\gamma,~ \Gamma$ both invertible, and the following equivalences according to Theorem \ref{thm:IsomeanV2CauchyMean}:
\begin{equation}\nonumber
    t=D_{lg}(a,b)=\isomeanvalue{x}{[a,~b]}{g,\gamma} ~~(l\in\int\gamma(x)g'(x)\mathrm{d}x),
\end{equation}
\begin{equation}\nonumber
    T=D_{LG}(a,b)=\isomeanvalue{x}{[a,~b]}{G,~\Gamma} ~~(L\in\int\Gamma(x)G'(x)\mathrm{d}x),
\end{equation}
($l,L$ are 3 times differentiable) and $f(D_{lg}(a,b))=\isomeanvalue{f(x)}{[a,~b]}{g,h}$, $f(D_{LG}(a,b))=\isomeanvalue{f(x)}{[a,~b]}{G,~H}$. Subsequently $l,g,L,G$ at least on $[a,b]$ meet the prerequisites of Losonczi's necessity Theorem \ref{thm:Los1}, according to which: (i)In the case $f$ is decreasing and $\isomeanvalue{f(x)}{[a,b]}{g,h}\leq \isomeanvalue{f(x)}{[a,b]}{G,H}$, which turns out $D_{lg}(a,b)\geq D_{LG}(a,b)$, it must hold the ``$\geq$'' case of (\ref{equ:IsomeanVNecess}) where $\gamma=l'/g'$ and $\Gamma=L'/G'$. (ii)In the case $f$ is increasing, analogously it must have ``$\leq$'' case of (\ref{equ:IsomeanVNecess}) for (\ref{equ:IsomeanC5Comp}).
}

\slrem{According to Losonczi's Remark 1, above theorem needs to assume $a,b$ are near enough; and if (\ref{equ:IsomeanVNecess}) is strict($<$ or $>$) then it's sufficient for (\ref{equ:IsomeanC5Comp}) to hold.}

Although Losonczi's necessity theorem need its generator functions $f,g$ be $(n+1)$ times differentiable, however the $g^{(n+1)}$ is prematurely cancellable item in its derivation, and does not affect the final condition's inequality. In our case only $l,L$ being 3 times differentiable are required.

In Section \ref{subsubsec:CompGeomElastic}, there is an example in which we work out how near $a,b$ to each other is necessary and sufficient for (\ref{equ:IsomeanC5Comp}).
\subsubsection{Sufficient conditions}
\slthm{\label{thm:IsomeanC5CompSuff}Let
(i) $f:[a,b]\to M$ be a strictly monotone, differentiable and bounded function;
(ii) $g,G$ be 2 times continuously differentiable bijections on $[a,b]$, on which $~g'\ne 0,~G'\ne 0$;
(iii) $h,H$ be 2 times continuously differentiable bijections on $M$, on which $~h'\ne 0,~H'\ne 0$. Then
\begin{equation} \label{equ:IsomeanVSufficient}
\begin{split}
\frac{\gamma(u)-\gamma(v)}{\gamma'(v)}~\frac{g'(u)}{g'(v)} \leq(\geq) \frac{\Gamma(u)-\Gamma(v)}{\Gamma'(v)}~\frac{G'(u)}{G'(v)}~~
\\ \big(u,v\in [a,b], ~~f~is~increasing(decreasing)\big)
\end{split}
\end{equation}
is sufficient for $\isomeanvalue{f(x)}{[a,b]}{g,h}\leq \isomeanvalue{f(x)}{[a,b]}{G,H}$(\ref{equ:IsomeanC5Comp}) to hold, where $\gamma:=h\circ f$, $\Gamma:=H\circ f$.
}
\slprf{With the same equivalences in previous proof and $\gamma=l'/g'$, $\Gamma=L'/G'$, the holding of (\ref{equ:IsomeanVSufficient})(``$\leq$'' case) is sufficient for $D_{lg}(a,b)\leq D_{LG}(a,b)$ according to Losonczi's sufficiency Theorem \ref{thm:Los2}. Therefore if with increasing $f$, $\isomeanvalue{f(x)}{[a,b]}{g,h}\leq \isomeanvalue{f(x)}{[a,b]}{G,H}$. The other case is analogous.
}

\subsubsection{Comparison of geometric mean and elastic mean}\label{subsubsec:CompGeomElastic}
Let $y=f(x)$ defined on $[a,b]~~(b>a>0)$ be positive, 2 times continuously differentiable and monotone. Here we denote the set of all $f$ as $CPM([a,b])$.

As an example, we shall find specific necessary conditions, and sufficient conditions for
\begin{equation}\label{equ:ElasGTGeom}
\isomeanvalue{f(x)}{[a,b]}{x,\ln y}\leq\isomeanvalue{f(x)}{[a,b]}{\ln x,y} ~~\big(f\in~CPM([a,b])\big).
\end{equation}
The left side is the geometric mean(will be denoted as $G$), and the right side is the elastic mean(denoted as $E$). (\ref{equ:ElasGTGeom}) is whereby denoted by $G\leq E$.

\paragraph{Necessary conditions\texorpdfstring{\\}{}}
Assume $f$ is increasing first. With (\ref{equ:IsomeanVNecess}) it is easy to check that:
\begin{equation}
\begin{split}
    \frac{(\ln f)''}{(\ln f)'} \leq \frac{f''}{f'}-\frac{2}{x}
\end{split}
\end{equation}
is necessary. It follows that $xf'/f\geq 2$ is necessary(note: $xf'/f$ is called the elasticity of $f$). If $f$ is decreasing then $xf'/f\leq 2$ is necessary.

Recalling the remark that if the inequalities are strict, then they are also sufficient for $\isomeanvalue{f(x)}{[a,b]}{g,h}\leq\isomeanvalue{f(x)}{[a,b]}{G,H}$ provided $a,b$ are near enough. Also in the case $f$ is decreasing, $xf'/f\leq 0 <2$, such we have the following propositions:

\slprop{Let $f\in CPM([a,b])$ and $a,b$ be provided near enough. If $f$ is increasing then its elasticity $xf'/f>2$ is necessary and sufficient for (\ref{equ:ElasGTGeom}). If $f$ is decreasing then (\ref{equ:ElasGTGeom}) shall hold.}

\slprop{For $f(x)=x^p$ on [a,b] with $a,b$ being near enough. If $p$ is positive, then $p>2$ is necessary and sufficient for (\ref{equ:ElasGTGeom}). If $p<0$ then (\ref{equ:ElasGTGeom}) holds.}


\paragraph{Sufficient conditions}
\sllem{$x^y<\exp{(x-1)}$ if 1 separates $x>0,~y>0$.}
\slprf{It is easy to prove that $\forall x(x\ne1),~\exp(x-1)>x$. (i)When $y>1>x>0$, $\exp(x-1)>x>x^y$. (ii)When $x>1>y>0$, also $\exp(x-1)>x>x^y$.}

\slthm{Let $f\in CPM([a,b])$. Then (\ref{equ:ElasGTGeom}) holds if
\begin{equation} \label{equ:ElasGTGeomSuff}
    \frac{f(u)}{f(v)}-\frac{u}{v}\ln\big(\frac{f(u)}{f(v)}\big)-1 \geq0, ~~(u,v\in [a,b]).
\end{equation}
}
\slprf{
Assume $f$ is increasing first. With (\ref{equ:IsomeanVSufficient}) it is easy to check that:
\begin{equation} \nonumber
\begin{split}
    \frac{(\ln f(u)-\ln f(v))f(v)}{f'(v)} \leq \frac{f(u)-f(v)}{f'(v)}~\frac{v}{u},~~\big(u,v\in [a,b])
\end{split}
\end{equation}
is sufficient. It follows that $(\ln f(u)-\ln f(v))f(v) \leq (f(u)-f(v))v/u$ is sufficient and obviously it's the same for decreasing $f$. Finally we get an unified (\ref{equ:ElasGTGeomSuff}) is sufficient regardless of $f$ being either increasing or decreasing.}

\slcor{\label{cor:DecFunGeoLessElas}Let $f\in CPM([a,b])$ be decreasing (e.g. $f(x)=x^p~(p<0)$), then (\ref{equ:ElasGTGeom}) holds.}
\slprf{
With (\ref{equ:ElasGTGeomSuff}), let $w=f(u)/f(v)$, $r=u/v$ then $w-r\ln w-1\geq0~\Rightarrow~w^r\leq\exp(w-1)$ is sufficient. Since $f$ is strictly decreasing, then $w,r$ are positive and separated by 1 or both are 1, such $w^r\leq\exp(w-1)$. The sufficient conditions are always met for decreasing $f$, no matter $a,b$ are near or far to each other.
}

\paragraph{Comprehensive comparison of \emph{E} and \emph{G} of \texorpdfstring{$x^p$}{power function}}
\sllem{\label{lem:funcPowerGra}For function $S(r)=r^p-p~r\ln r-1~(r>0,p\ne 0)$, a fixed (trivial) root for $S(r)=0$ is $r=1$, and\\
\indent1). If $p<0$, then $S(r)\ge0$ and 1 is the sole root;\\
\indent2). If $p>2$, there is sole second root $\alpha<1$ such that (a)$S(\alpha)=0$, (b)$S(r)<0$ ($r\in(0,\alpha)$), (c)$S(r)\geq0$ ($r\in[\alpha,+\infty)$);\\
\indent3). If $1<p<2$, there is sole second root $\beta>1$ such that (a)$S(\beta)=0$, (b)$S(r)\leq0$ ($r\in(0,\beta]$), (c)$S(r)>0$ ($r\in(\beta,+\infty)$);\\
\indent4). If $0<p\leq1$, then $S(r)\leq0$ and 1 is the sole root;\\
\indent5). If $p=2$ then (a)$S(r)<0$ ($r\in(0,1)$), (b)$S(r)>0$ ($r\in(1,+\infty)$).}

This is obtained by the analysis of $S$ via $S'$, the limit of $r^p/(p~r\ln r+1)$ etc. Then we have the following theorem.

\slthm{\label{thm:CompPowerEG}Let $f(x)=x^p~(x>0)$ defined on $[a,b]$ ($b>a>0$), and $\alpha<1,\beta>1$ be the roots of $S(r)$ in corresponding to case 2) and 3) resp. in Lemma \ref{lem:funcPowerGra}. Then listed are sufficient conditions for
\begin{equation}\label{equ:ElasGTGeom1}
\isomeanvalue{x^p}{[a,b]}{x,\ln y}\leq\isomeanvalue{x^p}{[a,b]}{\ln x,y} ~~~(G\leq E)
\end{equation}
or
\begin{equation}\label{equ:ElasLTGeom}
\isomeanvalue{x^p}{[a,b]}{\ln x,y}\leq\isomeanvalue{x^p}{[a,b]}{x,\ln y} ~~~(E\leq G)~:
\end{equation}
\indent1). If $p<0$, then (\ref{equ:ElasGTGeom1}) holds;\\
\indent2). If $p>2$ and $a<b\leq a/\alpha$ then (\ref{equ:ElasGTGeom1}) holds;\\
\indent3). If $1<p<2$ and $a<b\leq a\beta$ then (\ref{equ:ElasLTGeom}) holds;\\
\indent4). If $0<p\leq1$, then (\ref{equ:ElasLTGeom}) holds.
}
\slprf{Put $f$ in (\ref{equ:ElasGTGeomSuff}) and let $r=u/v$ we get: $G\leq E$ if the inequality
\begin{equation}
    S(r):=r^p-pr\ln r-1 \geq0, ~~(r\in[a/b,~b/a]).
\end{equation}
Another symmetrical theorem should be true, by which $E\leq G$ if
\begin{equation}
    S(r):=r^p-pr\ln r-1 \leq0, ~~(r\in[a/b,~b/a]).
\end{equation}
\indent1). In the case $p<0$ $f$ is decreasing, such $G\leq E$ is true with Corollary \ref{cor:DecFunGeoLessElas}. Also it's verified by Lemma \ref{lem:funcPowerGra} case 1) that $r^p-pr\ln r-1 \geq0$ if $p<0$.\\
\indent2). In the case $p>2$ and $\alpha<1$ being the root of $S(r)=0$ of Lemma \ref{lem:funcPowerGra} case 2), if $a<b\leq a/\alpha$ then $r\in[\alpha, 1/\alpha]$, on which $S(r)\geq 0$. It's sufficient for $G\leq E$.\\
\indent3). In the case $1<p<2$ and $\beta>1$ being the root of $S(r)=0$ of Lemma \ref{lem:funcPowerGra} case 3), if $a<b\leq a\beta$ then $r\in[1/\beta, \beta]$, on which $S(r)\leq 0$. It's sufficient for $E\leq G$.\\
\indent4). In the case $0<p\leq1$, according to Lemma \ref{lem:funcPowerGra} case 4) $S(r)\leq 0$. Such also $E\leq G$.
}

\slrem{All above are sufficient conditions, however not necessary ones. The above theorem is not yet a total solution for comparison of these 2 means. e.g. In case 2) $p>2$ we don't have a sufficient condition for otherwise $E\leq G$, since interval $[0,\alpha]$ does not contain any $r=u/v>1$, because $[0,\alpha]$ never spans over point 1. Similar are with case 3) and the case if $p=2$.}

\paragraph{Direct comparison of \emph{E} and \emph{G} of \texorpdfstring{$x^p$}{power function} by relative difference}
\slrem{By derivation of the formulae of geometric mean and elastic mean:
\begin{equation}
    G=\isomeanvalue{x^p}{[a,b]}{x,\ln y}=\exp\bigl(\big(\frac{b\ln b-a\ln a}{b-a}-1\big)p\bigl)=(I(a,b))^p,
\end{equation}
\begin{equation}
    E=\isomeanvalue{x^p}{[a,b]}{\ln x,y}=\frac{b^p-a^p}{p(\ln b-\ln a)}=(L_p(a,b))^p,
\end{equation}
where $I(a,b)$ is the identric mean of $a,b$ and $L_p(a,b)$ is the $p$-order logarithmic mean of $a,b$\cite{YANGZHCvCompHM}. Especially $p=2$, $G=(I(a,b))^2$ and $E=A(a,b)L(a,b)$ where $A(a,b)$ is the arithmetic mean and $L(a,b)$ is the logarithmic mean.}

\sldef{Let $y=x^p~(p\ne0)$ defined on $[a,b]$ ($b>a>0$). The number
\begin{equation}
\begin{split}
    \sigma_{GE}=(\isomeanvalue{x^p}{[a,b]}{x,\ln y}-\isomeanvalue{x^p}{[a,b]}{\ln x,y})/
                \isomeanvalue{x^p}{[a,b]}{\ln x,y}=G/E-1\\
\end{split}
\end{equation}
is called the relative difference of the geometric mean to the elastic mean of $x^p$.}

\sllem{With above definition if the ratio $r=b/a$ or $t=a/b$ is fixed, then \\
\indent1). $\sigma_{GE}$ is the function of $(r,p)$ or $(t,p)$:
\begin{equation}
\begin{split}
    \sigma_{GE}(r,p) =\frac{pr^{{\frac{pr}{r-1}}}\ln r}{e^p(r^p-1)}-1, ~~~
    \sigma_{GE}(t,p) =\frac{pt^{{\frac{pt}{t-1}}}\ln t}{e^p(t^p-1)}-1.
\end{split}
\end{equation}
\indent2). $\sigma_{GE}(r,p)=\sigma_{GE}(t,p)$. \\
\indent3). $\lim_{r\to1}\sigma_{GE}(r,p)=0$.
}
All above are easy to check and the proof is omitted.

\slthm{\label{thm:compEGbySigma}Let $y=x^p~(p\ne0)$ defined on $[a,b]$ ($b>a>0$) and $r=b/a$. The sufficient and necessary condition for $G<E$($E<G$) is $\sigma_{GE}(r,p)=\sigma_{GE}(1/r,p)<0(>0)$.}
\slprf{ By the definition and the lemma,
$\sigma_{GE}(r,p)=\sigma_{GE}(1/r,p)=(G-E)/E$, therefore $G<E \Leftrightarrow \sigma_{GE}(r,p)<0$ and $E<G \Leftrightarrow \sigma_{GE}(r,p)>0$ since $E>0$.}

\slcor{$\isomeanvalue{x^2}{[a,b]}{x,\ln y}>\isomeanvalue{x^2}{[a,b]}{\ln x,y}$~($G>E$) on $[a,b]$ ($b>a>0$).
}
\slprf{By checking up to 6 times derivatives of
\begin{equation}
\begin{split}
    \sigma_{GE}(r,2)=\frac{2r^{{\frac{2r}{r-1}}}\ln r}{e^2(r^2-1)}-1 
\end{split}
\end{equation}
 to $r$, or the derivatives of the factors of the derivatives, till it is no more transcendental, we are able to claim that $\sigma_{GE}(r,2)$ is decreasing when $r<1$ or is increasing when $r>1$, while $\lim_{r\to1}\sigma_{GE}(r,2)=0$. Such $\sigma_{GE}(r,2)>0 \Leftrightarrow G>E$.
}

With the theorem it's manually verifiable that Theorem \ref{thm:CompPowerEG} are sufficient conditions which are not sharp. For example, if $p=3$ then $\alpha=0.2142142...$ such if $a<b\leq a/0.2142142...$ then sufficiently $\isomeanvalue{x^3}{[a,b]}{x,\ln y}\leq\isomeanvalue{x^3}{[a,b]}{\ln x,y}$; however as far as $a<b\leq666a$ is still safe for the inequality, as with Theorem \ref{thm:compEGbySigma}, the threshold for the inequality to reverse near 666 can be manually verified by the graph of $\sigma_{GE}(r,3)$. This method is numerically efficient as the bivariate inequality turns to an univariate one($p$ is fixed).\\

The next subsections are all about the comparison methods for various scenarios, derived by help of monotonicity and convexity, and all are sufficient conditions.

\subsection{Comparison of 2 isomorphic means class I of a function}
The following theorem is derived from an integral extension form of Jensen's inequality. (Its proof is omitted.)
\linespread{1.1}
\slthm{\label{thm:IsoMeanC1Comp} Given $f\colon D\to M$ defined on interval $D$, $M\subseteq$ interval $I$, and 2  PVDMs $g,h$,
\\\indent1). If $g$ is increasing on $I$, and $g(h^{-1})$ is convex on $h(I)$, ~then $M_f|_g \geq M_f|_h$;
\\\indent2). If $g$ is increasing on $I$, and $g(h^{-1})$ is concave on $h(I)$, then $M_f|_g \leq M_f|_h$;
\\\indent3). If $g$ is decreasing on $I$, and $g(h^{-1})$ is convex on $h(I)$, ~then $M_f|_g \leq M_f|_h$;
\\\indent4). If $g$ is decreasing on $I$, and $g(h^{-1})$ is concave on $h(I)$, then $M_f|_g \geq M_f|_h$.
\\And with $g$, $h$ being derivable ($h'\ne0$),
\\\indent5). If $|g'/h'|$ is increasing on $I$, then $M_f|_g \geq M_f|_h$;
\\\indent6). If $|g'/h'|$ is decreasing on $I$, then $M_f|_g \leq M_f|_h$.
}\linespread{1.0}
The case 1) through 4) resembles Theorem \ref{thm:IsoMeanCompOld} (for numbers), also the Theorem 3 of \cite{LOSONLCAUCHYCOMP} which is about the case $g=G$ there. The case 5), 6) resembles Theorem \ref{thm:IsoMeanCompNew} (for numbers). However all cases are applicable to an ordinary $f$ on $D$, which may be discontinuous or non-monotone.

All 6 cases are sufficient but not necessary conditions, e.g. if with case 5) $D$ is extended just a little bit such it happens $|g'/h'|$ is not monotone on the extended $M$ (and $I$), it's very possible that still holds $M_f|_g \geq M_f|_h$ (on the new $D$), in which case the original $D$ may take the major weight.

\subsection{Comparison of 2 class II of a monotone function}
For comparison of class II, the monotonicity of the function $f$ is required in this section.

\linespread{1.1}
\slthm{\label{thm:IsoMeanC2CompCvx} Given a monotone $f$ defined on $[a,b]$ and 2 IVDMs $g,h$, in the cases $f$ increases,
\\\indent1). If $g$ increases on $[a,b]$, and $g(h^{-1})$ is convex on $h([a,b])$, ~then $M_f|_g^{II} \geq M_f|_h^{II}$;
\\\indent2). If $g$ increases on $[a,b]$, and $g(h^{-1})$ is concave on $h([a,b])$, then $M_f|_g^{II} \leq M_f|_h^{II}$;
\\\indent3). If $g$ decreases on $[a,b]$, and $g(h^{-1})$ is convex on $h([a,b])$, ~then $M_f|_g^{II} \leq M_f|_h^{II}$;
\\\indent4). If $g$ decreases on $[a,b]$, and $g(h^{-1})$ is concave on $h([a,b])$, then $M_f|_g^{II} \geq M_f|_h^{II}$.
\\In the cases $f$ decreases, all these 4 inequalities reverses. And with $g,h$ being differentiable ($h'\ne0$),
\\\indent5). If both of $|g'/h'|$ and $f$ increases or both decreases, then $M_f|_g^{II} \geq M_f|_h^{II}$;
\\\indent6). If one of $|g'/h'|$ and $f$ increases and another decreases, then $M_f|_g^{II} \leq M_f|_h^{II}$.
\\If $f$ is strictly monotone, all the inequalities in these 6 cases are strict($> or <$).
}
\linespread{1.2}
\slprf{Let
\begin{equation} \nonumber
\begin{split}
k =\frac{g(b)-g(a)}{h(b)-h(a)}\neq 0, ~
v &=h_2(x) =kh(x)-kh(a)+g(a), 
~u=h(x) ~~(x\in[a,b]).
\end{split}
\end{equation}
Such an $h_2\in\mathbb{V}h$ is said to be vertically aligned to $g(x)$ on $[a,b]$, since\\
\indent(i). $h_2(a)=g(a)$, $h_2(b)=g(b)$;\\
\indent(ii). interval $h_2([a,b])=g([a,b])$;\\
\indent(iii). $h_2,~g$ are of the same monotonicity;\\
\indent(iv). Based on above, $h_2^{-1}$ and $g^{-1}$ both ends at $(g(a),a), (g(b),b)$. $\Longrightarrow g(h_2^{-1}(v))$ ends at 2 points $(g(a),g(a)), (g(b),g(b))$ which are on the line $x=v$;\\
\indent(v). Besides, according to Lemma \ref{lem:HssVssMirror} $h_2\in\mathbb{V}h \Longrightarrow h_2^{-1}\in\mathbb{H}(h^{-1})$, therefore $g(h^{-1}(u))$ and $g(h_2^{-1}(v))$ are also H-scaleshifts, and have the same convexity according to Lemma \ref{lem:HssSameConvex}.\\
\indent Next we distinguish 4 cases to compare $f\circ h_2^{-1}(v)$ and $f\circ g^{-1}(v) ~(\forall v\in g([a,b]))$:\\
\indent1). $g(h^{-1}(u))$ is convex $\Longrightarrow g(h_2^{-1}(v))$ is convex on $g([a,b])\Longrightarrow g(h_2^{-1}(v))\leq v$; meanwhile $g$ is increasing and $f$ is increasing $\Longrightarrow f(h_2^{-1}(v)) \leq f(g^{-1}(v))$. According to Theorem \ref{thm:IsoMeanMonotone} (monotone property)\\
\begin{equation}
\begin{split}
&\isomeanvalue{f\circ h_2^{-1}}{[g(a),g(b)]}{\mathbb{V}v,\mathbb{V}y} \leq \isomeanvalue{f\circ g^{-1}}{[g(a),g(b)]}{\mathbb{V}v,\mathbb{V}y}\\
\Longrightarrow &M_f|_{h_2}^{II} \leq M_f|_g^{II} ~\overset{h_2\in\mathbb{V}h}{\Longrightarrow} ~M_f|_{h}^{II} \leq M_f|_g^{II}.
\end{split}
\end{equation}
\indent2). $g(h^{-1}(u))$ is concave $\Longrightarrow g(h_2^{-1}(v))$ is concave on $g([a,b])\Longrightarrow g(h_2^{-1}(v))\geq v$; meanwhile $g$ is increasing and $f$ is increasing $\Longrightarrow f(h_2^{-1}(v)) \geq f(g^{-1}(v))$. Similarly this leads to $M_f|_{h}^{II} \geq M_f|_g^{II}$.\\
\indent3). 4). ...(analogous and omitted.)\\
\indent So far we get 4 inequalities as the case 1) through 4). In the cases $f$ decreases while other conditions remain unchanged, it's obvious all these 4 inequalities reverses.\\
\indent As for case 5) and 6), $h_2$ is also differentiable. To examine the monotonicities of the derivative:
\begin{equation}
\begin{split}
\frac{\mathrm{d}g(h_2^{-1}(v))}{\mathrm{d}v} &= \frac{\mathrm{d}g(h_2^{-1}(v))}{\mathrm{d}h_2^{-1}(v)}\cdot \frac{\mathrm{d}h_2^{-1}(v)}{\mathrm{d}v} = g'\frac{1}{h_2'} = \frac{g'}{kh'}|_x = (\frac{g'}{kh'}\circ h_2^{-1})|_v
\end{split}
\end{equation}
in the following 4 cases and check their equivalences when combined with monotonicities of $h_2,g$(below $\nearrow$:~increasing, $\searrow$:~decreasing):\\
\indent1). When $(\frac{g'}{kh'}\circ h_2^{-1})|_v \nearrow$, $g(h_2^{-1}(v))\leq v$  $\overset{h_2\nearrow g\nearrow}{\Longleftrightarrow}$ when $\frac{g'}{kh'}\nearrow$,$h_2^{-1}(v)\leq g^{-1}(v)$;\\
\indent2). When $(\frac{g'}{kh'}\circ h_2^{-1})|_v \nearrow$, $g(h_2^{-1}(v))\leq v$  $\overset{h_2\searrow g\searrow}{\Longleftrightarrow}$ when $\frac{g'}{kh'}\searrow$,$h_2^{-1}(v)\geq g^{-1}(v)$;\\
\indent3). When $(\frac{g'}{kh'}\circ h_2^{-1})|_v \searrow$, $g(h_2^{-1}(v))\geq v$  $\overset{h_2\nearrow g\nearrow}{\Longleftrightarrow}$ when $\frac{g'}{kh'}\searrow$,$h_2^{-1}(v)\geq g^{-1}(v)$;\\
\indent4). When $(\frac{g'}{kh'}\circ h_2^{-1})|_v \searrow$, $g(h_2^{-1}(v))\geq v$  $\overset{h_2\searrow g\searrow}{\Longleftrightarrow}$ when $\frac{g'}{kh'}\nearrow$,$h_2^{-1}(v)\leq g^{-1}(v)$.\\
From there 4 cases can merge to 2 cases without conflictions. And knowing the monotonicity of $g'/(kh')$ is the same as that of $|g'/h'|$, these summarize as:\\
\indent a). When $|g'/h'|$ increases, $h_2^{-1}(v)\leq g^{-1}(v)$;\\
\indent b). When $|g'/h'|$ decreases, $h_2^{-1}(v)\geq g^{-1}(v)$.\\
By further combining $f$ and applying Theorem \ref{thm:IsoMeanMonotone}, case 5) and 6) are proved.\\
\indent Finally if $f$ is strictly monotone then $f\circ h_2^{-1}, f\circ g^{-1}$ are strictly monotone, and based on convexity of $g(h_2^{-1}(v))$ for every case there is at least 1 sub-interval of g([a,b]) such that on which $f(h_2^{-1}(v)) >(or <) f(g^{-1}(v))$ while for the rest parts $f(h_2^{-1}(v)) = f(g^{-1}(v))$~(Otherwise the contradiction is that there is no convexity of $g(h_2^{-1}(v))$ at all). By applying the strict inequality of Theorem \ref{thm:IsoMeanMonotone} or its corollary, all the inequalities are strict. \\
\indent This completes the proof of Theorem \ref{thm:IsoMeanC2CompCvx}.
}
\linespread{1.0}
\slrem{For easier practice, let $m=|k|/k=sgn(k)$, then $|g'/h'|$=$mg'/h'$.}

There is a trivial out-of-scope case when $mg'/h'$ is a constant $C>0$. Then by antiderivative we have $mg=Ch+C'$ hence two isomorphic means are equal. Also noticeable is that the monotonicity of $|g'/h'|$ is invariant with V-scaleshifts of $g$, $h$.

There are some examples of comparison of class II:
\begin{itemize} \setlength{\itemsep}{-0.1em}
  \item[(1).] $f(x)=\tan x, ~g(x)=\ln x, ~h(x)=x, ~x\in[0.1, 1.5]$. $|g'/h'|=\ln'x/x'=1/x$, It is decreasing. $\tan x$ is increasing. Thus according to Theorem \ref{thm:IsoMeanC2CompCvx} case 6), $\isomeanvalueII{\tan x}{[0.1, 1.5]}{\ln x} \leq \isomeanvalueII{\tan x}{[0.1, 1.5]}{x}=\overline {{\displaystyle {\tan x}}_{\scriptstyle {[0.1,1.5]}}}$. Also have noted that $\isomeanvalueII{\tan x}{(0, ~\pi/2)}{\ln x} =2/\pi$, while $\overline {{\displaystyle {\tan x}}_{\scriptstyle {(0, ~\pi/2)}}}$ does not exist.
  \item[(2).] $f(x)=x, ~g(x)=x^2, ~h(x)=x, ~x\in[a, b], ~b>a\geq0$. $|g'/h'|=(x^2)'/x'=2x$, It is increasing. $f$ is increasing. Thus according to Theorem \ref{thm:IsoMeanC2CompCvx} case 5), $\isomeanvalueII{x}{[a,b]}{x^2} \geq \isomeanvalueII{x}{[a,b]}{x}=\overline {{\displaystyle {x}}_{\scriptstyle {[a,b]}}}$. Simple calculation leads to: ~For positive $a,~b$,
  \begin{equation}
    \frac{2(a^2+ab+b^2)}{3(a+b)} \geq \frac{a+b}{2}.
  \end{equation}
  \item[(3).] $f(x)=x, ~g(x)=\sin x, ~h(x)=\cos x, ~x\in[a,b], ~0\leq a<b\leq\pi/2$. $|g'/h'|=-\sin'x/\cos'x=\cot x$, It is decreasing. $f$ is increasing. Thus according to Theorem \ref{thm:IsoMeanC2CompCvx} case 6), $\isomeanvalueII{x}{[a,b]}{\sin x} \leq \isomeanvalueII{x}{[a,b]}{\cos x}$. After derivations, including those after $\overline {{\displaystyle {x}}_{\scriptstyle {[a,b]}}}$ participating in comparisons, we get: ~For $0\leq a\leq\pi/2$, $~0\leq b\leq\pi/2(~b\neq a)$,
  \begin{equation}
    \frac{b\sin b-a\sin a+\cos b-\cos a}{\sin b-\sin a} < \frac{a+b}{2}
    <~\frac{a\cos a-b\cos b+\sin b-\sin a}{\cos a-\cos b}.
  \end{equation}
\end{itemize}

\subsection{Comparison of 2 class IV of a monotone \emph{f} -- ~a partial solution}
\slthm{\label{thm:IsoMeanCompC4p}
Let $g,G$ be 2 differentiable IVDMs and $h,H$ be 2 differentiable PVDMs of monotone $f\colon [a,b]\to M$. $G'\ne0$ on $[a,b]$, $H'\ne0$ on $I=[\min\{M\},\max\{M\}]$.
\setlength{\parskip}{-0.4em}
\begin{itemize}\setlength{\itemsep}{-0.2em}
  \item [1).] If both $|g'/G'|$ and $f$ are increasing or decreasing on $[a,b]$, and $|h'/H'|$ is increasing on $I$, then $M_f|_{g,h} \geq M_f|_{G,H}$.
  \item [2).] If one of $|g'/G'|$ and $f$ is increasing another decreasing on $[a,b]$, and $|h'/H'|$ is decreasing on $I$, then $M_f|_{g,h} \leq M_f|_{G,H}$.
\end{itemize}
If $f$ is strictly monotone, these 2 inequalities are strict($> or <$).
}

\slprf{With the aid of vertically alignment of a V-scaleshift $G_2$ of $G$ to $g$ on $[a,b]$ same as that in the proof of previous theorem, the proof can be done mainly via cooperation of Theorem \ref{thm:IsoMeanMonotone} and Theorem \ref{thm:IsoMeanC1Comp}, since $M_f|_{g,h}$ and $M_f|_{G_2,h}=M_f|_{G,h}$ can be treated as 2 isomorphic means class I on the same $g([a,b])$. For e.g. case 1). it's therefore easy to prove $M_f|_{g,h} \geq M_f|_{G,h} \geq M_f|_{G,H}$.}

 The detailed proof should initially list 4 cases, while only these 2 cases therein are of decidable comparison results. Hence it is only a partial solution as there are 2 other cases not being able to handled by the theorem. The methods of proving in the next 5 subsections are similar and most proofs are omitted.

\subsection{Comparison of 2 class V -- ~a partial solution}
As a special case of previous problem, this is to compare $M_x|_{g,h}, ~M_x|_{G,H}$: while $f(x)=x$ being strictly increasing, the inequalities are strict in below corollary.
\slcor{Let $g,G$ be 2 differentiable IVDMs and $h,H$ be 2 differentiable PVDMs of $f(x)=x~(x\in[a,b])$. $G'\ne0$, $H'\ne0$ on $[a,b]$.
\setlength{\parskip}{-0.4em}
\begin{itemize}\setlength{\itemsep}{-0.2em}
  \item [1).] If both $|g'/G'|$ and $|h'/H'|$ are increasing, then $M_x|_{g,h} > M_x|_{G,H}$.
  \item [2).] If both $|g'/G'|$ and $|h'/H'|$ are decreasing, then $M_x|_{g,h} < M_x|_{G,H}$.
\end{itemize}
}

\subsection{Comparison of 2 class IV with only different IVDMs\texorpdfstring{\\($g\ne G,h=H)$}{}}
This is to compare $M_f|_{g,h},~M_f|_{G,h}$. Like comparing isomorphic means class II, $f$'s monotone is required in below theorem.
\slthm{Let $g,G$ be 2 differentiable IVDMs and $h$ be a PVDM of monotone $f$ defined on $[a,b]$. $G'\ne0$ on $[a,b]$.
\setlength{\parskip}{-0.4em}
\begin{itemize}\setlength{\itemsep}{-0.2em}
  \item[1).] If both $|g'/G'|$ and $f$ increases or both decreases, then $M_f|_{g,h} \geq M_f|_{G,h}$;
  \item[2).] If one of $|g'/G'|$ and $f$ increases another decreases, then $M_f|_{g,h} \leq M_f|_{G,h}$.
\end{itemize}
If $f$ is strictly monotone, these 2 inequalities are strict($> or <$).
}

\subsection{Comparison of 2 class IV with only different PVDMs \texorpdfstring{\\($g=G,h\ne H$)}{}}
This is to compare $M_f|_{g,h},~M_f|_{g,H}$, as if comparing 2 isomorphic means class I. $f$'s monotone is NOT required in below theorem.
\slthm{Let $f\colon D\to M$ be defined on interval $D$, $M\subseteq$ interval $I$. Let $g$ be an IVDM and $h,H$ be 2 differentiable PVDMs of $f$. $H'\ne0$ on $I$.
\setlength{\parskip}{-0.4em}
\begin{itemize}\setlength{\itemsep}{-0.2em}
  \item [1).] If $|h'/H'|$ is increasing on $I$, then $M_f|_{g,h} \geq M_f|_{g,H}$.
  \item [2).] If $|h'/H'|$ is decreasing on $I$, then $M_f|_{g,h} \leq M_f|_{g,H}$.
\end{itemize}
}

\subsection{Comparison of 2 class III of an increasing \emph{f}}
\slthm{\label{thm:IsoCompC3Inc} Let 2 differentiable mappings $g,h$ be both IVDMs and PVDMs of an increasing $f\colon [a,b]\to M$. $h'\ne0$ on $[a,b]$ and $I=[\min\{M\},\max\{M\}]$.
\setlength{\parskip}{-0.4em}
\begin{itemize}\setlength{\itemsep}{-0.2em}
  \item [1).] If $|g'/h'|$ is increasing on $[a,b]$ and $I$ then $M_f|_{g,g} \geq M_f|_{h,h}$.
  \item [2).] If $|g'/h'|$ is decreasing on $[a,b]$ and $I$ then $M_f|_{g,g} \leq M_f|_{h,h}$.
\end{itemize}
If $f$ is strictly monotone, these 2 inequalities are strict($> or <$).
}
\slprf{Based on proof of Theorem \ref{thm:IsoMeanCompC4p}, 4 cases could be listed initially with $f$ possibly being either increasing or decreasing. Then 2 cases with decreasing $f$ are of undecidable comparison results. The rest 2 cases make the theorem true.}

The simplest application of above is to compare 2 isomorphic means class III of $f(x)=x$ on $[a,b]$ generated respectively by $g(x)=x^p~(p\ne0)$ and by $h(x)=x^q~(q\ne0)$. If $p>q$ then $|g'/h'|=|p/q|x^{(p-q)}$ is increasing, therefore $M_f|_{g,g} > M_f|_{h,h} \Longrightarrow \big(\frac{a^p+b^p}{2}\big)^{(1/p)} > \big(\frac{a^q+b^q}{2}\big)^{(1/q)}$ for $~a\ne b, p>q$.

\subsection{Comparison of 2 class IV of a decreasing \emph{f} with exchanged DMs}
\slthm{Let 2 differentiable mappings $g,h$ be both IVDMs and PVDMs of a decreasing $f\colon [a,b]\to M$. $h'\ne0$ on $[a,b]$ and $I=[\min\{M\},\max\{M\}]$.
\setlength{\parskip}{-0.4em}
\begin{itemize}\setlength{\itemsep}{-0.2em}
  \item [1).] If $|g'/h'|$ is increasing on $[a,b]$ and $I$ then $M_f|_{g,h} \leq M_f|_{h,g}$.
  \item [2).] If $|g'/h'|$ is decreasing on $[a,b]$ and $I$ then $M_f|_{g,h} \geq M_f|_{h,g}$.
\end{itemize}
If $f$ is strictly monotone, these 2 inequalities are strict($> or <$).
}
\slprf{Based on proof of Theorem \ref{thm:IsoMeanCompC4p}, 4 cases could be listed initially with $f$ possibly being either increasing or decreasing. Then 2 cases with increasing $f$ are of undecidable comparison results. The rest 2 cases make the theorem true.}
Below are 3 examples:
\begin{enumerate}\setlength{\itemsep}{-0.2em}
  \item[(1).] By choosing $g(x)=\cos x$, $h(x)=\sin x$, $f(x)=\pi/2-x$, $[a,b]\subseteq[0,\pi/2]$, and applying this method, we finally conclude with a bivariate inequality, to both sides of which further applying the inverse function of $f$ yields
\begin{equation} 
    \arccos(\frac{\cos a+\cos b}2) >  \arcsin(\frac{\sin a+\sin b}2).
\end{equation}
  \item[(2).] By choosing $g(x)=x^p$,~$p\ne0$, $h(x)=x^q$,~$q\ne0$, $f(x)=1/x$, $[a,b]\subseteq(0,+\infty)$ and applying the theorem, we get an inequality, to both sides of which further applying $f^{-1}$ yields
\begin{equation} 
    \biggl( \frac{(p-q)(b^p-a^p)}{p(b^{p-q}-a^{p-q})}\biggl)^{1/q} <
    \biggl( \frac{(q-p)(b^q-a^q)}{q(b^{q-p}-a^{q-p})}\biggl)^{1/p} ~for ~p<q.
\end{equation}
The above reverses for $p>q$.
  \item[(3).] By choosing $g(x)=\ln x$, $h(x)=x$, and a decreasing $f(x)>0$ on $[a,b]\subseteq(0,+\infty)$, we get a result similar to Corollary \ref{cor:DecFunGeoLessElas}, i.e. for a decreasing positive function (has not to be continuous), $G\leq E$.
\end{enumerate}

\pagebreak

\section{Conclusion \& vision}
As a specific topic related to functions bonded on the isomorphic frames, isomorphic means are uncovered as a huge family of mean values. Various types of means of numbers and of function can be treated as the special cases, instances or derivations of isomorphic means. Such it crosses with some existing concepts of means, e.g. geometric mean of a function, Stolarsky means, Cauchy mean values and its derived forms.

The isomorphic means can be used for derivation of inequalities, especially when the comparison problems of its various sub-classes can be solved in a systematical way.

The root concept: the isomorphic frames and its basic derivations: the isomorphic number and the DVI function can be deemed as the fundamentals of some type of ``mathematics related to isomorphic frames'', which covers isomorphic means. The ``dual-variable-isomorphic convex function'' as introduced in reference article \cite{LIUY} (which is also based on DVI function though it never refined the isomorphic frames) can also be deemed as an example of such mathematics, and is closely related to the isomorphic means. Just for extra information, the ``geometrically convex function'' is also a special case of it, as discussed in \cite{LIUY}. Therefore the vision of this article is to expand the scope of MA in this similar way by tying of more MA concepts to isomorphic frames and studying their extended properties and behaviors, e.g. the ``dual-variable-isomorphic convexity''.

To facilitate it by precisely graphically representing the isomorphic frames and its embedded and bonded objects, the so-called ``$n$-dimensional isomorphic coordinate system'' could have been introduced as another root concept, which generally represents a type of uneven space rendered by the isomorphic frame. The devise of such coordinate system is supported by the Lemma \ref{lem:OrderSetbeBJ}, which validates a bijection that maps the space onto a part of a Cartesian coordinate system. However this is abridged from this article for the time being.

On the coordinate systems the geometrical meaning of isomorphic means and of ``dual-variable-isomorphic convex function'' can be demonstrated and correlated. One also can see that the geometrical meaning of Cauchy mean value in the coordinate system is different than that of isomorphic means.  On the systems, more isomorphic-frame-related mathematics can be introduced, e.g. the so-called densities on the coordinate systems, the graphs of functions, the convexity of sets and functions, the slopes, lengths, curvatures, even the so-called local densities of the graphs of functions on the systems can be discussed.

The author is looking forward to presenting more theories of extension of this genre, which are all expected to establish on the basis of $n$-dimensional isomorphic frames. Meaningfully for current scope, let's end by trying to give a premature definition of isomorphic mean of a function of ($n-1$) variables($n\ge2$):

\sldef{Let intervals $X_1,...,X_n,U_1,...,U_n\subseteq\mathbb{R}$, and $g_i\colon X_i\to U_i(i=1,...,n)$ be $n$ continuous and monotone bijections. Function $f:D\to M$ of ($n-1$) variables is bounded, and $f\wedge\mathscr{I}_m\{g_1,...,g_n\}$. $D$ and $E=\mathscr{I}_m\{g_1,...,g_{n-1}\}(D)$ each are connected sets and are measurable of $(n-1)$ dimensional hypervolume. ~If there exists $M_\varphi \in U_n$, being the mean value of $\varphi\colon=g_n\circ f(g_1^{-1},...,{g_{n-1}}^{-1})\colon E\to U_n$ on $E$, ~then $g_n^{-1}(M_\varphi)\in X_n$ ~is called the (all-variable-) isomorphic mean of $f$ on $D$ generated by $\mathscr{I}_m\{g_1,...,g_n\}$, denoted by $\isomeanvalue{f}{D}{g_1,...,g_n}$,
\begin{equation}\label{equ:IsoMeanValueN}
    \isomeanvalue{f}{D}{g_1,...,g_n} =
    g_n^{-1}\biggl(\frac{\int_{\scriptscriptstyle E}g_n\circ f(g_1^{-1},...,{g_{n-1}}^{-1})\mathrm{d}^{(n-1)}u}
    {\int_{\scriptscriptstyle E}\mathrm{d}^{(n-1)}u}\biggl). \\
\end{equation}
}

\pagebreak

\end{document}